%% file: Main_File-Stochastic_Sampling.tex
\documentclass[review, 12pt]{elsarticle}
\usepackage{lineno,hyperref}
\modulolinenumbers[5]
\usepackage{setspace}
\usepackage[justification=centering]{caption}
\usepackage{subfigure}
\usepackage{amssymb}
\usepackage{color}
\usepackage{amsmath}
\usepackage{nomencl}
\usepackage{array}
\usepackage{multirow}
\usepackage{geometry}
\usepackage{makeidx}
\usepackage{float}
\DeclareMathOperator{\E}{\mathbb{E}}
\DeclareMathOperator{\Trace}{trace}
\DeclareMathOperator{\Var}{Var}
\DeclareMathOperator{\Dev}{Dev}
\makeindex
\makeindex
\providecommand{\printnomenclature}{\printglossary}
\providecommand{\makenomenclature}{\makeglossary}

\usepackage{etoolbox}
\makenomenclature

\def\bm{\boldsymbol}

\begin{document}

\begin{frontmatter}

\title{Stochastic Sampling for Structural Topology Optimization with Many Load Cases: Density-Based and Ground Structure Approaches}

\author[a]{Xiaojia (Shelly) Zhang}
\ead{xzhang645@gatech.edu}

\author[b]{Eric de Sturler}
\ead{sturler@vt.edu}

\author[a]{Glaucio H. Paulino}
\ead{paulino@gatech.edu}

\address[a]{School of Civil and Environmental Engineering, Georgia Institute of Technology, 790 Atlantic Drive, Atlanta, GA, 30332, USA}
\address[b]{Department of Mathematics, Virginia Tech, McBryde Hall, 225 Stanger Street, Blacksburg, VA, 24061, USA}

\begin{abstract}

We propose an efficient probabilistic method to solve a deterministic problem -- we present a randomized optimization approach that drastically reduces the enormous computational cost of optimizing designs under many load cases for both continuum and truss topology optimization. Practical structural designs by topology optimization typically involve many load cases, possibly hundreds or more. The optimal design minimizes a, possibly weighted, average of the compliance under each load case (or some other objective). This means that in each optimization step a large finite element problem must be solved for each load case, leading to an enormous computational effort. On the contrary, the proposed randomized optimization method with stochastic sampling requires the solution of only a few (e.g., 5 or 6) finite element problems (large linear systems) per optimization step. Based on simulated annealing, we introduce a damping scheme for the randomized approach. Through numerical examples in two and three dimensions, we demonstrate that the stochastic algorithm drastically reduces computational cost to obtain similar final topologies and results (e.g., compliance) compared with the standard algorithms. The results indicate that the damping scheme is effective and leads to rapid convergence of the proposed algorithm.

\end{abstract}


\end{frontmatter}
\input{Section1-Introduction_v1}

\input{Section2-StochasticSampling_v1}

\input{Section3-AlgorithmicParameters_v1}

\input{Section4-NumericalExamples_v1}

\input{Section5-ConcludingRemarks_v1}

\section*{Acknowledgements}
The authors G. H. Paulino and X. Zhang acknowledge the financial support from the US National Science Foundation (NSF) under projects \#1559594 (formerly \#1335160) and \#1321661. We are also grateful for the endowment provided by the Raymond Allen Jones Chair at the Georgia Institute of Technology. The work by E. de Sturler was supported in part by the grant NSF DMS 1217156. The information provided in this paper is the sole opinion of the authors and does not necessarily reflect the view of the sponsoring agencies.

\appendix
\nomenclature{$n_f$}{Number of monitored step for the discrete filter in the stochastic algorithm}
\nomenclature{$n_{\text{step}}$}{Sample window size in the damping scheme}
\nomenclature{$\tau_{\text{step}}$}{Tolerance for the damping scheme}
\nomenclature{$\gamma$}{Scale factor in the damping scheme}
\nomenclature{$R$}{Effective step ratio in the damping scheme}
\nomenclature{$\alpha_\text{Top}$}{Resolution of the structural topology}
\nomenclature{$\alpha_i$}{Weighting factor for the $i$th load case}
\nomenclature{$\alpha_f$}{Discrete filter value}
\nomenclature{$\boldsymbol{f}_i$}{External force vector the $i$th load case}
\nomenclature{$\boldsymbol{u}_i$}{Displacement vector the $i$th load case}
\nomenclature{$\boldsymbol{\xi}$}{Random vector with its entries drawn from the Rademacher distribution }
\nomenclature{$\boldsymbol{\rho}^{*}$}{Optimal solution obtained by the standard algorithm in the density-based method}
\nomenclature{$\boldsymbol{x}^{*}$}{Optimal solution obtained by the standard algorithm in the GSM}
\nomenclature{$\hat{\boldsymbol{\rho}}^{S}$}{Optimal solution obtained by the stochastic algorithm in the density-based method}
\nomenclature{$\hat{\boldsymbol{x}}^{S}$}{Optimal solution obtained by the stochastic algorithm in the GSM}
\nomenclature{$N_{\text{solve}}$}{Total number of linear systems solves in the otpimization process}
\nomenclature{$\bm{H}$}{Filter matrix for density}
\nomenclature{$\overline{\boldsymbol{\rho}}$}{Filtered density vector}
\nomenclature{$p$}{Penalization parameter in the density-based method}
\nomenclature{$d$}{Number of dimensions}
\nomenclature{$E_0$}{Young's modulus of the solid material}
\nomenclature{$\textbf{\textit{F}}$}{Weighted external force matrix}
\nomenclature{$C$}{Weighted compliance (objective function)}
\nomenclature{$\hat{C}$}{Estimated objective function by the stochastic algorithm}
\nomenclature{$\textbf{\textit{K}}$}{Global stiffness matrix}
\nomenclature{$\rho^{(e)}$}{Density of element $e$ ($e$th design variable in the density-based method)}
\nomenclature{$L^{(e)}$}{Length of truss member $e$}
\nomenclature{$M$}{Number of elements in the mesh/initial ground structure}
\nomenclature{$m$}{Number of load cases}
\nomenclature{$N$}{Number of nodes in the mesh/initial ground structure}
\nomenclature{$n$}{Sample size}
\nomenclature{$\tau_{\text{opt}}$}{Tolerance for the optimization process}
\nomenclature{$\textbf{\textit{U}}$}{Weighted displacement matrix}
\nomenclature{$v^{(e)}$}{Volume of element $e$}
\nomenclature{$V_\text{max}$}{Prescribed maximum volume}
\nomenclature{$x^{(e)}$}{Cross-sectional area of member $e$ ($e$th design variable in the GSM)}\
\nomenclature{$x^{\text{min}}, x^{\text{max}}$}{Lower and upper bounds for the cross-sectional area of the member in the GSM}
\nomenclature{$\rho^{\text{min}}$}{Lower bound for density in the density-based method}
\section{Nomenclature}
\vspace*{-0.3in}
\printnomenclature

\newpage

\bibliographystyle{elsarticle-num}
\bibliography{Shelly_ref}

\end{document}

%% file: Section1-Introduction_v1.tex
\section{Introduction}
Structural topology optimization is an important tool that allows for improved designs with minimal design iterations. In the field of structural topology optimization, designs accounting for many load cases are practical and provide stable structures. Indeed, real-world structural designs, for example, high-rise buildings and long-span bridges, generally involve numerous load cases \cite{diaz1992,bendsoe1994,bental1997,bendsoe2003}. For the end-compliance minimization formulation, several methods that include many load cases have been published \cite{bendsoe2003}. The main idea is to minimize a proper norm of the vector of compliances from all load cases. This paper concentrates on one popular method, the weighted sum formulation that minimizes a weighted sum of the vector of compliances from all load cases, which requires the solution of the structural equation for each load case. One advantage of this formulation is that the weighting corresponds to the probability of the load occurring \cite{christiansen2001}, which naturally connects to the field of uncertainty quantification. Another formulation that accounts for many load cases is the min-max formulation. This formulation minimizes the worst-case compliance, or the infinity-norm, of all load cases \cite{achtziger1998,bendsoe2003}. This formulation also requires the solutions of linear equations from all the load cases, therefore, in terms of computational cost, it is similar to the aforementioned weighted sum formulation. In this paper, we proposed a randomized algorithm that drastically reduces the computational cost. Similar techniques have been applied to solve inverse problems -- see \cite{haber2012} and references within, and in the context of the solutions of least squares problems -- see \cite{avron2010} and references within.   
  
The remainder of the paper is organized as follows. The remainder of Section 1 describes the weighted sum formulation of topology optimization with many load cases. We briefly review the weighted sum formulation of the standard nested minimum end-compliance topology optimization with many load cases, using both the density based method and the ground structure method. Section 2 reviews the stochastic sampling techniques we use to estimate the trace of a general matrix and proposes the stochastic algorithm for minimum end-compliance topology optimization under many load cases. Section 3 discusses the algorithmic parameters and introduces a damping scheme for the stochastic algorithm. Section 4 presents numerical examples in two- and three- dimensions highlighting the efficiency and effectiveness of the proposed algorithm, and Section 5 provides concluding remarks with suggestions for extending the work.

\subsection{Density-based topology optimization formulation with many load cases}


For a set of $m$ given load cases $\boldsymbol{f}_{i}, i=1,...,m$, the standard weighted sum formulation for the minimum end-compliance design with many load cases using the density-based method can be written as \cite{bendsoe2003},
\begin{equation}\label{eq:density_obj}
\begin{aligned}
\min_{\boldsymbol{\rho}} C\left(\boldsymbol{\rho}\right)&=\min_{\boldsymbol{\rho}}\sum_{i=1}^{m}\alpha_i\bm{f}_i^{\,\,T}\bm{u}_i(\boldsymbol{\rho}) , \\
\text{s.t.}\quad & \sum_{e=1}^M v^{(e)}\bar{\rho}^{(e)}-V_{\text{max}}\leq0 \,, \\
 & 0<\rho^{\text{min}}\leq\rho^{(e)}\leq1,\,e=1,...,M , \\
 & \text{with } {\bm{u}}_i(\bar{\boldsymbol{\rho}}) = \textit{\textbf{K}}\left(\boldsymbol{E}\left(\bar{\boldsymbol{\rho}}\right)\right)^{-1}\textit{\textbf{f}}_i,\, i=1,...,m\,.\end{aligned}
\end{equation}
In this minimization problem, the objective function $C$ is the weighted-average compliance of the corresponding structure, $\boldsymbol{\rho}\in\mathbb{R}^{M}$ is the vector of design variables (the density field), and $M$ is the number of elements in the finite element mesh. We define $\bar{\boldsymbol{\rho}}$ and $\boldsymbol{H}$ as the filtered physical density and the filter matrix such that $\bar{\boldsymbol{\rho}}=\boldsymbol{H}\boldsymbol{\rho}$ \cite{Bourdin2001}. In order to ensure a positive definite stiffness matrix $\textit{\textbf{K}}\in\mathbb{R}^{dN\times dN}$, a lower bound $\rho^{\text{min}}$ is prescribed on $\rho^{(e)}$, where $N$ is the number of nodes and $d$ is the dimension of the problem, so $dN$ is the number of degrees of freedom. The volume (area) of element $e$ is given by $v^{(e)}$, and $V_{\text{max}}$ is the prescribed upper bound on the total volume. The weight and the displacement vector associated with load case $\textit{\textbf{f}}_i\in\mathbb{R}^{dN}$ are given by $\alpha_i$ ($\alpha_i>0,\, \sum_{i=1}^{m}\alpha_i=1$) and $\textit{\textbf{u}}_i\in\mathbb{R}^{dN}$, respectively. The Young's modulus $\boldsymbol{E}$ is defined by, for example, the Solid Isotropic Material with Penalization (SIMP) \cite{bensoe1989,zhou1991} approach. Other models, e.g., RAMP (Rational Approximation of Material Properties) \cite{bendsoe2003,stolpe2001}, can be used and would not alter the conceptual presentation of the topic. For the SIMP approach with the density filter, we have $E^{(e)}=E_\text{min}+[\bar{\rho}^{(e)}]^p (E_0-E_\text{min})$, where $E_0$ and $E_\text{min}$ are elastic moduli for solid material and Ersatz material, respectively, and $p$ is a penalization factor. 
The gradient (sensitivity) of the objective function corresponds to a weighted sum of the sensitivities of each individual loading cases, which can be expressed as
\begin{equation}\label{eq:density_sens}
 \nabla C_\rho^{(e)}=\frac{\partial C}{\partial \rho^{(e)}}=-\sum_{i=1}^{m}\alpha_i\bm{u}_i^{\,\,T}\frac{\partial \textit{\textbf{K}}}{\partial \rho^{(e)}}\textit{\textbf{u}}_i.
\end{equation}
The formulation \eqref{eq:density_obj} is known to be convex when $p=1$ \cite{petersson1999}. Using $p > 1$ to obtain a solid-void solution, one makes the problem non-convex and, as expected, the solution obtained from the optimization algorithm may not be the global minimum.

\subsection{Ground-structure based formulation with many load cases}
The standard weighted sum formulation for the minimum end-compliance design with the ground structure method that accounts for many load cases takes the following form,
\begin{equation}\label{eq:ground_obj}
\begin{aligned}
\min_{\boldsymbol{x}} C\left(\boldsymbol{x}\right)&=\min_{\boldsymbol{x}}\sum_{i=1}^{m}\alpha_i\bm{f}_i^{\,\,T}\textit{\textbf{u}}_i(\boldsymbol{x}) ,\\
\text{s.t.}\quad & \sum_{e=1}^M L^{(e)} x^{(e)}-V_{\text{max}}\leq0 \,,\\
 & 0<x^{\text{min}}\leq x^{(e)}\leq x^{\text{max}},\,e=1,...,M ,\\
 & \text{with } \bm{u}_i(\boldsymbol{x})=\bm{K}\left(\boldsymbol{x}\right)^{-1}\bm{f}_i,\, i=1,...,m\,.\end{aligned}
\end{equation}
The vector $\boldsymbol{x}\in\mathbb{R}^{M}$ is a vector of design variables, with component $x^{(e)}$ being the cross-sectional area of truss member $e$. It is subject to lower bound $x^{\text{min}}$ and upper bound $x^{\text{max}}$. Furthermore, $M$ is the number of truss members in the ground structure, $L^{(e)}$ is the length of truss member $e$, and $V_{\text{max}}$ is the prescribed upper bound on the total volume. For the $i$th load case $\textit{\textbf{f}}_i\in\mathbb{R}^{dN}$, $\alpha_i$ and $\textit{\textbf{u}}_i\in\mathbb{R}^{dN}$ are the corresponding weight factor and displacement vector. As in the density-based method, the sensitivity of the objective function in the ground structure method is the weighted sum of the sensitivity from each load case:
\begin{equation}\label{eq:ground_sens}
 \nabla C_x^{(e)}=\frac{\partial C}{\partial x^{(e)}}=-\sum_{i=1}^{m}\alpha_i\bm{u}_i^{\,\,T}\frac{\partial \textit{\textbf{K}}}{\partial x^{(e)}}\bm{u}_i.
\end{equation}
The standard nested formulation of the end-compliance objective function with a single load case has been proven to be convex in \cite{Svanberg1984} for a positive definite stiffness matrix and in \cite{Achtziger1997} for a positive semi-definite stiffness matrix. The formulation \eqref{eq:ground_obj} with multiple load cases is easily proved to be convex.
 
\subsection{Synopsis} 
The optimization process contains three main components: solving the structural equilibrium problem for a set of given design variables, computing the gradient of the objective function, and updating the design variables. In this paper, we use the density-based and ground structure methods, both methods utilize the Optimality Criterion (OC) algorithm \cite{groenwold2008} as the update scheme, which only requires gradient information of the objective function and the volume constraint.
The convergence criterion for optimization used in the formulations \eqref{eq:density_obj} and \eqref{eq:ground_obj} is that the maximum change in design variables drops below a given tolerance, $\tau_\text{opt}$.

The overall computational cost of the standard optimization formulations \eqref{eq:density_obj} and \eqref{eq:ground_obj} can be estimated by the total number of structural equations (large linear system) solves, i.e., $m\times N_\text{step}$, where $N_\text{step}$ is the number of optimization steps. For realistic three-dimensional (3D) problems, where the mesh size, the number of design variables, and the number of load cases are large, the associated computational cost is enormous. In this paper, we propose a randomized approach that reduces the problem with many load cases ($m$) to a sequence of optimization steps with only a few load cases ($n \ll m$) to solve per step, that yields results similar to those from the standard method.

%


%% file: Section2-StochasticSampling_v1.tex
\section{Stochastic sampling and topology optimization}
This section proposes a stochastic sampling approach for the minimum end-compliance topology optimization formulations with many load cases. First, we briefly review a stochastic sampling technique to estimate the trace of a general matrix. Furthermore, the stochastic sampling technique is applied to the minimum end-compliance topology optimization with many load cases using both density-based and ground structure methods.
 
\subsection{Stochastic sampling of matrices}
For a general matrix $\textit{\textbf{A}}\in\mathbb{R}^{q\times q}$,  stochastic sampling techniques can be used to estimate the trace of $\textit{\textbf{A}}$. We discuss one popular approach here, the Hutchinson trace estimator \cite{hutchinson1989}, but several alternatives exist, e.g., the Gaussian estimators and unit vector estimators -- see \cite{avron2011,roosta2015} and references therein. Let $\boldsymbol{\xi}\in\mathbb{R}^q$ be a random vector containing entries that are independent and identically distributed (i.i.d.) with value $\pm1$ each with probability $1/2$. This distribution is known as the Rademacher distribution.
It follows immediately that, for each entry, the expectation  $\E\left(\xi_i\right)=0$. Since the entries are independent, the expectation of $\xi_i\xi_j$ is given by,
\begin{equation} 
\E\left(\xi_i\xi_j\right)=\left\{\begin{array}{ll}
               {0\,,}  & {\rm i\neq j\,,}\\
							\\
               {1\,,} & {\rm i=j\,.}
             \end{array}
             \right.
\end{equation} 
Now consider the expectation of the random variable $\boldsymbol{\xi}^T\bm{A}\boldsymbol{\xi}$,
\begin{equation}
\E\left(\boldsymbol{\xi}^T\bm{A}\boldsymbol{\xi}\right)=\E\left(\sum_{i=1}^{q}\sum_{i=1}^{q}\xi_i A_{ij}\xi_j\right)=\sum_{i=1}^{q}\sum_{i=1}^{q}A_{ij}\E\left(\xi_i\xi_j\right).
\end{equation}
Since $\E\left(\xi_i\xi_j\right)=1$ only when $i=j$ and 0 otherwise, we get
\begin{equation}
\E\left(\boldsymbol{\xi}^{T}\bm{A}\boldsymbol{\xi}\right)=\sum_{i=1}^{q}A_{ii}=\Trace\left(\textit{\textbf{A}}\right)\,.
\end{equation}
As a result, the trace of a given matrix $\textit{\textbf{A}}$ can be estimated using random samples. Here, we utilize the sample average approximation (SAA) technique, see, e.g., \cite{shapiro2009}, which approximates the expected value by the average. The sample average or empirical mean for $n$ samples (or realizations) $\boldsymbol{\xi}_1, \boldsymbol{\xi}_2, ..., \boldsymbol{\xi}_n$, of the random variable $\boldsymbol{\xi}$ is defined as
\begin{equation}
 \E_{S}\left(\boldsymbol{\xi}^T\textit{\textbf{A}}\boldsymbol{\xi}\right)=\frac{1}{n}\sum_{k=1}^n\boldsymbol{\xi}_k^T\textit{\textbf{A}}\boldsymbol{\xi}_k\,.
\end{equation}
According to the Law of Large Numbers (LLN), as the number of (independent) samples $n$ approaches $\infty$, the sample average $\E_S\left(\boldsymbol{\xi}^T\textit{\textbf{A}}\boldsymbol{\xi}\right)$ converges to the expectation $\E\left(\boldsymbol{\xi}^T\textit{\textbf{A}}\boldsymbol{\xi}\right)$, which is the trace of $\bm{A}$\,,
%
\begin{equation}
\E_{S}\left(\boldsymbol{\xi}^T\textit{\textbf{A}}\boldsymbol{\xi}\right)=\frac{1}{n}\sum_{k=1}^n\boldsymbol{\xi}_k^T\textit{\textbf{A}}\boldsymbol{\xi}_k \rightarrow \E\left(\boldsymbol{\xi}^{T}\textit{\textbf{A}}\boldsymbol{\xi}\right)=\Trace\left(\textit{\textbf{A}}\right)\,.
\end{equation}
The use of trace estimators can be highly efficient if the matrix \textit{\textbf{A}} is not explicitly available and its computation is expensive. We also have that $\E_S\left(\boldsymbol{\xi}^T\textit{\textbf{A}}\boldsymbol{\xi}\right)$ is an unbiased estimator of $\E\left(\boldsymbol{\xi}^T\textit{\textbf{A}}\boldsymbol{\xi}\right)$ \cite{shapiro2009}. Since the expected accuracy of such estimates is given by the variance, the goal is to obtain the estimates with the smallest variance. The Rademacher distribution, from which we choose the random vectors $\boldsymbol{\xi}$, was shown in \cite{hutchinson1989} to be the distribution that yields the smallest variance. Other studies rank various trace estimators differently according to criteria other than the variance.
For further information, readers are referred to, e.g., \cite{avron2011,roosta2015}. For a symmetric matrix \textit{\textbf{A}}, the variance and standard deviation of $\boldsymbol{\xi}^T\textit{\textbf{A}}\boldsymbol{\xi}$ are defined as
\begin{equation}
\Var\left(\boldsymbol{\xi}^T\textit{\textbf{A}}\boldsymbol{\xi}\right)=\E\left\{\left[\boldsymbol{\xi}^T\textit{\textbf{A}}\boldsymbol{\xi}-\E\left(\boldsymbol{\xi}^T\textit{\textbf{A}}\boldsymbol{\xi}\right)\right]^2\right\}=2\sum^{q}_{\substack{{i,j=1}\\ {i\neq j}}}A_{ij}^2\,,
\end{equation}
and
\begin{equation}
\Dev\left(\boldsymbol{\xi}^T\textit{\textbf{A}}\boldsymbol{\xi}\right)=\sqrt{\Var\left(\boldsymbol{\xi}^T\textit{\textbf{A}}\boldsymbol{\xi}\right)}=\sqrt{2\sum_{\substack{{i,j=1}\\ {i\neq j}}}^{q}A_{ij}^2}\,.
\end{equation}
%
The variance and standard deviation can be estimated using a finite number of samples. Given the samples above, we define the sample variance and sample standard deviation as
\begin{equation}
\Var_S\left(\boldsymbol{\xi}^T\textit{\textbf{A}}\boldsymbol{\xi}\right)=\frac{1}{n}\sum_{k=1}^{n}\left[\boldsymbol{\xi}_k^T\textit{\textbf{A}}\boldsymbol{\xi}_k-\frac{1}{n}\sum_{\ell=1}^{n}\boldsymbol{\xi}_\ell^T\textit{\textbf{A}}\boldsymbol{\xi}_\ell\right]^2\,,
\end{equation}
and
\begin{equation}
\Dev_S\left(\boldsymbol{\xi}^T\textit{\textbf{A}}\boldsymbol{\xi}\right)=\sqrt{\frac{1}{n}\sum_{k=1}^{n}\left[\boldsymbol{\xi}_k^T\textit{\textbf{A}}\boldsymbol{\xi}_k-\frac{1}{n}\sum_{\ell=1}^{n}\boldsymbol{\xi}_\ell^T\textit{\textbf{A}}\boldsymbol{\xi}_\ell\right]^2}\,.
\end{equation}
Below, we use these derivations for estimating the objective function (the compliance), as well as its gradient or sensitivity. A key issue for efficiency is that both can be estimated with the same set of random vectors.

\subsection{Randomized topology optimization with stochastic sampling}
We apply the stochastic sampling technique to both the density-based method and
the ground structure method. The basic idea is to replace the compliance and its
gradient by stochastic estimates and to use these estimates in the optimization algorithm.

\subsubsection{Density-based method}
Consider the standard topology optimization formulation in \eqref{eq:density_obj} with $m$ load cases 
$\bm{f}_i$,
$i=1,...,m$, and corresponding weights $\alpha_i$. We define a weighted load matrix
$\bm{F} \in \mathbb{R}^{dN\times m}$ as
$\bm{F} = \left[\sqrt{\alpha_1}\, \bm{f}_1,...,\sqrt{\alpha_m}\,\ \bm{f}_m\right]$.
In a similar fashion, we define the weighted displacement matrix
$\bm{U} \in \mathbb{R}^{dN\times m}=\left[\sqrt{\alpha_1}\,\, \textit{\textbf{u}}_1,...,\sqrt{\alpha_m}\,\,\textit{\textbf{u}}_m\right]$,
whose columns are the corresponding displacement fields.
The matrix $\textit{\textbf{U}}$ is defined by the equilibrium equation, $\textit{\textbf{U}}=\textit{\textbf{K}}^{\,\,-1}\textit{\textbf{F}}$.
So, we can write the end-compliance and its sensitivities as
traces of the symmetric matrices
\begin{equation}
C\left(\boldsymbol{\rho}\right)=\sum_{i=1}^{m}\alpha_i\textit{\textbf{f}}_i^{\,\,T}\textit{\textbf{u}}_i=\Trace\left(\textit{\textbf{F}}^{\,T}\textit{\textbf{U}}\right)=\Trace\left(\textit{\textbf{F}}^{\,T}\textit{\textbf{K}}^{\,\,-1}\textit{\textbf{F}}\right)
\end{equation}
and
\begin{equation}
 \nabla C_{\boldsymbol{\rho }}^{(e)}=-\sum_{i=1}^{m}\alpha_i\textit{\textbf{u}}^{\,\,T}_i\frac{\partial \textit{\textbf{K}}}{\partial \rho^{(e)}}\textit{\textbf{u}}_i=-\Trace\left(\textit{\textbf{U}}^{\,\,T}\frac{\partial \textit{\textbf{K}}}{\partial \rho^{(e)}}\textit{\textbf{U}}\right)=-\Trace\left(\textit{\textbf{F}}^{\,T}\textit{\textbf{K}}^{-1}\frac{\partial \textit{\textbf{K}}}{\partial \rho^{(e)}}\textit{\textbf{K}}^{-1}\textit{\textbf{F}}\right).
\end{equation}

With the random variable $\boldsymbol{\xi}$ as defined above,
the end-compliance and the topology optimization problem can be expressed as
\begin{equation}\label{eq:density_com}
C\left(\boldsymbol{\rho}\right) = \Trace\left(\textit{\textbf{F}}^{\,T}\textit{\textbf{K}}^{-1}\textit{\textbf{F}}\right)=\E\left(\boldsymbol{\xi}^T\textit{\textbf{F}}^{\,T}\textit{\textbf{K}}^{-1}\textit{\textbf{F}}\boldsymbol{\xi}\right)=\E\left[\left(\textit{\textbf{F}}\boldsymbol{\xi}\right)^T\textit{\textbf{K}}^{-1}\left(\textit{\textbf{F}}\boldsymbol{\xi}\right)\right],
\end{equation}
and
\begin{equation}\label{eq:density_formulation_s}
	\begin{aligned}
		\min_{\boldsymbol{\rho}} C\left(\boldsymbol{\rho}\right)&=\min_{\boldsymbol{\rho}} \E\left[\left(\textit{\textbf{F}}\boldsymbol{\xi}\right)^T\textit{\textbf{K}}\left(\boldsymbol{\rho}\right)^{-1}\left(\textit{\textbf{F}}\boldsymbol{\xi}\right)\right] ,\\
		\text{s.t.}\quad & \sum_{e=1}^M v^{(e)}\rho^{(e)}-V_{\text{max}}\leq0 \,,\\
		& 0<\rho^{\text{min}}\leq\rho^{(e)}\leq1,\,e=1,...,M\,.\\
	\end{aligned}
\end{equation}
The sensitivity of the objective function can be expressed as
\begin{equation}\label{eq:density_true_sens}
\nabla C_ {\boldsymbol{\rho }} ^{(e)} = -\Trace\left(\textit{\textbf{F}}^{\,T}\textit{\textbf{K}}^{-1}\frac{\partial \textit{\textbf{K}}}{\partial \rho^{(e)}}\textit{\textbf{K}}^{-1}\textit{\textbf{F}}\right)
=-\E\left[\left(\textit{\textbf{F}}\boldsymbol{\xi}\right)^{\,T}\textit{\textbf{K}}^{-1}\frac{\partial \textit{\textbf{K}}}{\partial \rho^{(e)}}\textit{\textbf{K}}^{-1}\left(\textit{\textbf{F}}\boldsymbol{\xi}\right)\right].
\end{equation}
Equations \eqref{eq:density_formulation_s} and \eqref{eq:density_true_sens}, stated as a stochastic programming
problem, are equivalent to the standard formulation \eqref{eq:density_obj} and \eqref{eq:density_sens}.
We approximate the compliance and its gradient by replacing their expectation in
the equations above by its sample average estimate. Given the i.i.d. random sample
$\boldsymbol{\xi}_1, \boldsymbol{\xi}_2, ..., \boldsymbol{\xi}_n$ as $n$ realizations of the random vector $\boldsymbol{\xi}$,
we define
\begin{equation}\label{eq:density_obj_sto}
 \widehat{C}^S\left(\boldsymbol{\rho}\right)=\frac{1}{n}\sum_{k}^{n}\left(\textit{\textbf{F}}\boldsymbol{\xi}_k\right)^T\textit{\textbf{K}}\left(\boldsymbol{\rho}\right)^{-1}\left(\textit{\textbf{F}}\boldsymbol{\xi}_k\right),
\end{equation}
and
\begin{equation}\label{eq:density_sen_sto}
(\nabla \widehat{C}_{\boldsymbol{\rho }}^{S})^{(e)} = \frac{\partial \widehat{C}^S}{\partial \rho^{(e)}}=\frac{1}{n}\sum_{k=1}^{n}\left(\textit{\textbf{F}}\boldsymbol{\xi}_k\right)^T\textit{\textbf{K}}^{-1}\frac{\partial \textit{\textbf{K}}}{\partial \rho^{(e)}}\textit{\textbf{K}}^{-1}\left(\textit{\textbf{F}}\boldsymbol{\xi}_k\right).
\end{equation}
We remark that $\widehat{C}^S\left(\boldsymbol{\rho}\right)$ and $\nabla \widehat{C}_{\boldsymbol{\rho }}^{S}$
are unbiased estimators for the compliance and its gradient.
By the LLN, $\widehat{C}^{S}\left(\boldsymbol{\rho}\right) \rightarrow C\left(\boldsymbol{\rho}\right)$ and
$\nabla \widehat{C}_{\boldsymbol{\rho }}^{S} \rightarrow \nabla C_{\boldsymbol{\rho }}$
(with probability 1)
for any feasible $\boldsymbol{\rho}$ as $n\rightarrow\infty$ \cite{anton2001}.
%

In each step of the optimization algorithm, we use the same random vectors to
estimate the compliance and its gradient using \eqref{eq:density_obj_sto} and \eqref{eq:density_sen_sto}.
To avoid convergence for a specific random load case, a new set of $n$ random vectors is selected
at each optimization step.
If $n \ll m$, our proposed algorithm reduces the computational cost from $m\times N_\text{step}$ to roughly $n\times N_\text{step}$ if the convergence of the optimization is not affected.

The idea to approximate the gradient and the objective function in a
structural optimization problem is similar to the use of
stochastic gradient-based methods in other areas, e.g., stochastic gradient
descent (SGD) \cite{robbins1951}, and stochastic meta-descent (SMD) \cite{schraudolph1999},
which are optimization methods mainly for unconstrained optimization problems.
Stochastic gradient methods use small sub-samples (also referred as mini-batches)
to estimate the gradient and have been applied in other fields,
such as large-scale machine learning \cite{vishwanathan2006,schraudolph2007}.
For our application, the SAA gradient estimator
always has a positive angle with the gradient,
and hence is a (negative) descent direction for the unconstrained case.
We demonstrate numerically how effective the approximation is in Section 4.
Because the optimization problem for the density-based method with penalization is non-convex,
deterministic gradient-based methods may not converge to the global minimum (as expected).
The stochastic algorithm leads to solutions that are (roughly) as accurate in terms of the end-compliance
as those from the standard algorithm.
This is demonstrated by the results in Section 4.2.


In general, for stochastic gradient methods, a damping or averaging scheme
(also referred to as a decay schedule of the scalar gain or gain vector adaptation)
is needed to achieve convergence \cite{schraudolph2007}.
Various damping or averaging methods for step sizes have been proposed for different types of problems \cite{nemirovski2009}.
For the structural optimization problems in this paper, we propose a damping scheme that adjusts the move limits (reminiscent of a trust region), which especially befits the structural optimization framework.
The idea of the proposed damping scheme is similar to adjusting the size of the update in the simulated annealing \cite{kirkpatrick1983,salamon2002} and is discussed in detail in Section \ref{dampingscheme}.
Robbins and Monro \cite{robbins1951} have given conditions for the convergence of stochastic
methods that use such damping schemes; see also \cite{schraudolph2007}.

The results in Section 4 show that the convergence of the proposed algorithm is typically rapid if our damping scheme is properly used, roughly as fast as for the standard algorithm and sometimes faster.
Since the randomized algorithm solves only $n$ linear systems at each optimization step,
compared with $m$ for the standard algorithm, and $n \ll m$,
then the proposed stochastic algorithm is computationally much more efficient.
Moreover, randomized algorithms with a proper damping scheme
can be more robust in finding the global minimum for non-convex optimization problems than deterministic algorithms.
We demonstrate this with numerical examples in Section 4.
In general, such increased robustness comes at the price of slower convergence.

To analyze the effects of randomization on the
stochastic optimization algorithms, we consider the variance and sample variance
for the compliance estimate and its gradient.
%
\begin{equation}\label{eq:den_var_obj}
\mathrm{Var}
  \left[ (\bm{F}\bm{\xi})^T \bm{K}^{-1} (\bm{F}\bm{\xi}) \right] = 2\sum^{m}_{\substack{{i,j=1}\\ {i\neq j}}}
  \left( \bm{F}^T \bm{K}^{-1} \bm{F} \right)^2_{ij}\,,
\end{equation}
and
\begin{equation}\label{eq:den_var_sen}
\mathrm{Var} \left[\left(\textit{\textbf{F}}\boldsymbol{\xi}\right)^T\textit{\textbf{K}}^{-1}\frac{\partial \textit{\textbf{K}}}{\partial \rho^{(e)}}\textit{\textbf{K}}^{-1}\left(\textit{\textbf{F}}\boldsymbol{\xi}\right)\right] = 2 \sum^{m}_{\substack{{i,j=1}\\ {i\neq j}}}\left( \bm{F}^T \bm{K}^{-1} \frac{\partial \bm{K}}{\partial \rho^{(e)}}\bm{K}^{-1} \bm{F} \right)^2_{ij}.
\end{equation}
Similarly, the sample variance of the compliance and its gradient can be expressed as
\begin{equation}\label{eq:den_sam_var_obj}
\mathrm{Var}_S\left[\left(\textit{\textbf{F}}\boldsymbol{\xi}\right)^T\textit{\textbf{K}}^{-1}\left(\textit{\textbf{F}}\boldsymbol{\xi}\right)\right]=\frac{1}{n}\sum_{k=1}^{n}\left(\boldsymbol{\xi}_k^T\textit{\textbf{F}}^T\textit{\textbf{K}}^{-1}\textit{\textbf{F}}\boldsymbol{\xi}_k-\frac{1}{n}\sum_{\ell=1}^{n}\boldsymbol{\xi}_\ell^T\textit{\textbf{F}}^T\textit{\textbf{K}}^{-1}\textit{\textbf{F}}\boldsymbol{\xi}_\ell\right)\,,
\end{equation}
and
\begin{multline}\label{eq:den_sam_var_sen}
\hspace*{0.5in}
\mathrm{Var}_S\left[\left(\textit{\textbf{F}}\boldsymbol{\xi}\right)^T\textit{\textbf{K}}^{-1}\frac{\partial \textit{\textbf{K}}}{\partial \rho^{(e)}}\textit{\textbf{K}}^{-1}\left(\textit{\textbf{F}}\boldsymbol{\xi}\right)\right]\\=\frac{1}{n}\sum_{k=1}^{n}\left(\boldsymbol{\xi}_k^T\textit{\textbf{F}}^{\,T}\textit{\textbf{K}}^{-1}\frac{\partial \textit{\textbf{K}}}{\partial \rho^{(e)}}\textit{\textbf{K}}^{-1}\textit{\textbf{F}}\boldsymbol{\xi}_k-\frac{1}{n}\sum_{\ell=1}^{n}\boldsymbol{\xi}_\ell^T\textit{\textbf{F}}^{\,T}\textit{\textbf{K}}^{-1}\frac{\partial \textit{\textbf{K}}}{\partial \rho^{(e)}}\textit{\textbf{K}}^{-1}\textit{\textbf{F}}\boldsymbol{\xi}_\ell\right)\,.
\end{multline}

\subsubsection{Ground structure method}
A stochastic version of the optimization problem for the ground structure method with the weighted force matrix $\bm{F}$ is analogous to that for the density-based method \eqref{eq:density_formulation_s}, which is given by
\begin{equation}\label{eq:ground_formulation_s}
\begin{aligned}
\min_{\boldsymbol{x}} C\left(\boldsymbol{x}\right)&=\min_{\boldsymbol{x}} \E\left[\left(\textit{\textbf{F}}\boldsymbol{\xi}\right)^T\textit{\textbf{K}}\left(\boldsymbol{x}\right)^{-1}\left(\textit{\textbf{F}}\boldsymbol{\xi}\right)\right] ,\\
\text{s.t.}\quad & \sum_{e=1}^M L^{(e)}x^{(e)}-V_{\text{max}}\leq0 \,,\\
& 0< x^{\text{min}}\leq x^{(e)}\leq x^{\text{max}},\,e=1,...,M\,,\\
\end{aligned}
\end{equation}
and the sensitivity of the objective function can be expressed as,
\begin{equation}\label{eq:ground_true_sens}
\nabla C_{\boldsymbol{x}}^{(e)}=-\E\left[\left(\textit{\textbf{F}}\boldsymbol{\xi}\right)^{\,T}\textit{\textbf{K}}^{-1}\frac{\partial \textit{\textbf{K}}}{\partial x^{(e)}}\textit{\textbf{K}}^{-1}\left(\textit{\textbf{F}}\boldsymbol{\xi}\right)\right].
\end{equation}
%

Here too, we replace the objective function and its gradient by their sample average estimate.
Using the same assumptions, for an i.i.d. random sample  $\boldsymbol{\xi}_1, \boldsymbol{\xi}_2, ..., \boldsymbol{\xi}_n$, the compliance $C\left(\boldsymbol{x}\right)$ and its gradient can be estimated by
\begin{equation}\label{eq:ground_obj_sto}
 \widehat{C}^S\left(\boldsymbol{x}\right)=\frac{1}{n}\sum_{k}^{n}\left(\textit{\textbf{F}}\boldsymbol{\xi}_k\right)^T\textit{\textbf{K}}\left(\boldsymbol{x}\right)^{-1}\left(\textit{\textbf{F}}\boldsymbol{\xi}_k\right)\,,
\end{equation}
and
\begin{equation}
(\nabla \widehat{C}_{\boldsymbol{x}}^{S})^{(e)}=\frac{\partial \widehat{C}^S}{\partial x^{(e)}}=-\frac{1}{n}\sum_{k=1}^{n}\left(\textit{\textbf{F}}\boldsymbol{\xi}_k\right)^T\textit{\textbf{K}}^{-1}\frac{\partial \textit{\textbf{K}}}{\partial x^{(e)}}\textit{\textbf{K}}^{-1}\left(\textit{\textbf{F}}\boldsymbol{\xi}_k\right)\,,
\end{equation}
where $\widehat{C}^S\left(\boldsymbol{x}\right)$ and $(\nabla \widehat{C}_{\boldsymbol{x}}^{S})^{(e)}$ are the estimated compliance and estimated sensitivity of the corresponding structure subjected to the sample load cases, which are analogous to \eqref{eq:density_obj_sto} and \eqref{eq:density_sen_sto}. The variance and sample variance of the objective function and its sensitivity take the same form as \eqref{eq:den_var_obj}-\eqref{eq:den_sam_var_sen}, and therefore, are not listed here for the sake of conciseness. Note that since the standard optimization formulation with the GM is convex, its optimal solution is a global minimum, hence the optimal values of all stochastic cases are equal to or larger than the one obtained from the standard formulation. This is confirmed by the results in Section 4.

We conclude this subsection by making a remark on the flexibility of our proposed randomized algorithm for topology optimization. This algorithm can be combined with many types of update schemes that are based on gradient information, such as the OC \cite{groenwold2008} and the Method of Moving Asymptotes (MMA) \cite{svanberg1987}.  In this work, we adopt OC as the update scheme for both continuum and truss topology optimization.


\subsection{Discrete filter for ground structure method with stochastic sampling}\label{section:DiscreteFilter} 
  
In the GSM, the discrete filter proposed in \cite{Ramos2016} is applied to the truss topology optimization with the stochastic algorithm to extract valid structures out of ground structures. The use of the filter reduces the number of redundant bars, which reduces the cost of subsequent optimization steps and helps to limit the effects of the stochastic estimates for the optimization problem. 
For the standard multiple load case optimization problem \eqref{eq:ground_obj}, the discrete filter can be expressed as
\begin{equation}\label{eq:orignal_filter}
Filter\left(\boldsymbol{x}_k,\alpha_{f}\right)=\left\{\begin{array}{ll}
               {0,}  & {\textbf{if}\quad\frac{\textit{x}^{(e)}_k}{\max\left(\boldsymbol{x}_k\right)}<\alpha_f<1\,,}\\
							\\
               {x^{(e)},} & {\rm otherwise\,,}
             \end{array}
             \right.
\end{equation}
where $\alpha_f$ is the prescribed filter value and $x_k^{(e)}$ is the design variable for truss member $e$ at the $k$th step. The discrete filter is applied to the ground structure at each step and removes the truss members with areas smaller than the filter $\alpha_f$.

For the stochastic case of the GSM, we form a discrete filter that slightly differs from the standard discrete filter. To limit the effects of stochastic estimates in the randomized algorithm, the discrete filter removes truss members only when their areas have remained below the prescribed filter value for $n_f$ steps, namely, 
\begin{equation}\label{eq:stoch_filter}
Filter^S\left(\boldsymbol{x}_k,\alpha_{f}\right)=\left\{\begin{array}{ll}
               {0,}  & {\text{if}\quad \max\left(\frac{\textit{x}^{(e)}_{k-n_f}}{\max\left(\boldsymbol{x}_{k-n_f}\right)},...,\frac{\textit{x}^{(e)}_{k}}{\max\left(\boldsymbol{x}_{k}\right)}\right)<\alpha_f<1},\\
							\\
               {x^{(e)},} & {\rm otherwise,}
             \end{array}
             \right.
\end{equation}
where $Filter^S(\boldsymbol{x}_k,\alpha_f)$ denotes the randomized filter (i.e., the discrete filter for the ground structure method with stochastic sampling), and $n_f$ is a chosen number of monitored steps (see discussions of parameters below). The use of the randomized filter ($Filter^S$) leads to a more efficient ground structure method with stochastic sampling than if the filter were not used (see Numerical Examples section).

%% file: Section3-AlgorithmicParameters_v1.tex
\section{A damping scheme and algorithmic parameters for randomized optimization}
This section introduces a damping scheme to achieve convergence and demonstrates the effectiveness of this scheme through a three-truss example. We further discuss the algorithmic parameters that are used in the proposed randomized optimization framework and comment on the range of values chosen for those parameters.
  
\subsection{The proposed damping scheme: effective step ratio and step size reduction}\label{dampingscheme}
In the randomized algorithm, the structure must be adjusted based on the random linear combination of load cases that changes at each optimization step. Therefore, the convergence criteria commonly used for the standard structural optimization framework is insufficient for the proposed framework. Based on the simulated annealing \cite{kirkpatrick1983,salamon2002}, we propose a damping scheme that evaluates the progress of steps and reduces the move limit (similar to the scalar gain in SGD) accordingly throughout the proposed randomized optimization framework. We simply define the effective step ratio $R$ as follows (using the GSM notation):
%
\begin{equation}\label{eq:stepev}
	R=\frac{\frac{1}{n_\text{step}}||\left(\boldsymbol{x}_k-\boldsymbol{x}_{k-n_\text{step}+1}\right)||}{||\boldsymbol{x}_k-\boldsymbol{x}_{k-1}||},
\end{equation}
where $n_\text{step}$ is the sample window size. The average step size over a sample window is divided by the current step length. 
This effective step ratio serves as an indicator of the optimizer's status, i.e., the ratio is relatively large when the optimizer is making progress; and relatively small (typically smaller than 0.1) when the step is not effective. The latter may be caused by two situations: the moving average (nominator) is small, then the optimizer is not making progress overall; or the current step length (denominator) is large, and this may be because the current step is dominated by extreme load cases. Once $R$ is below a prescribed tolerance $\tau_\text{step}$, i.e., $R<\tau_\text{step}$, we reduce the allowable move limit of the optimizer by a prescribed scale factor $\gamma$ to alleviate the influence from the random load cases. We start to reduce the move limit after $n_\text{step}$ steps.

The effectiveness of the damping scheme is illustrated through a simple numerical example for the GSM. This example is selected on purpose to show the poor performance of the stochastic optimization algorithm without a damping scheme. We consider a three-truss structure supported at their left ends and subjected to a set of 9 equal-weighted load cases, $\textit{\textbf{f}}_1,...,\textit{\textbf{f}}_9$, at their right joint, as shown in Fig. \ref{Fig1}. The problem formulation is given as follows:
%
\begin{figure}[!htbp]
	\centering
	\includegraphics[width=6.2in]{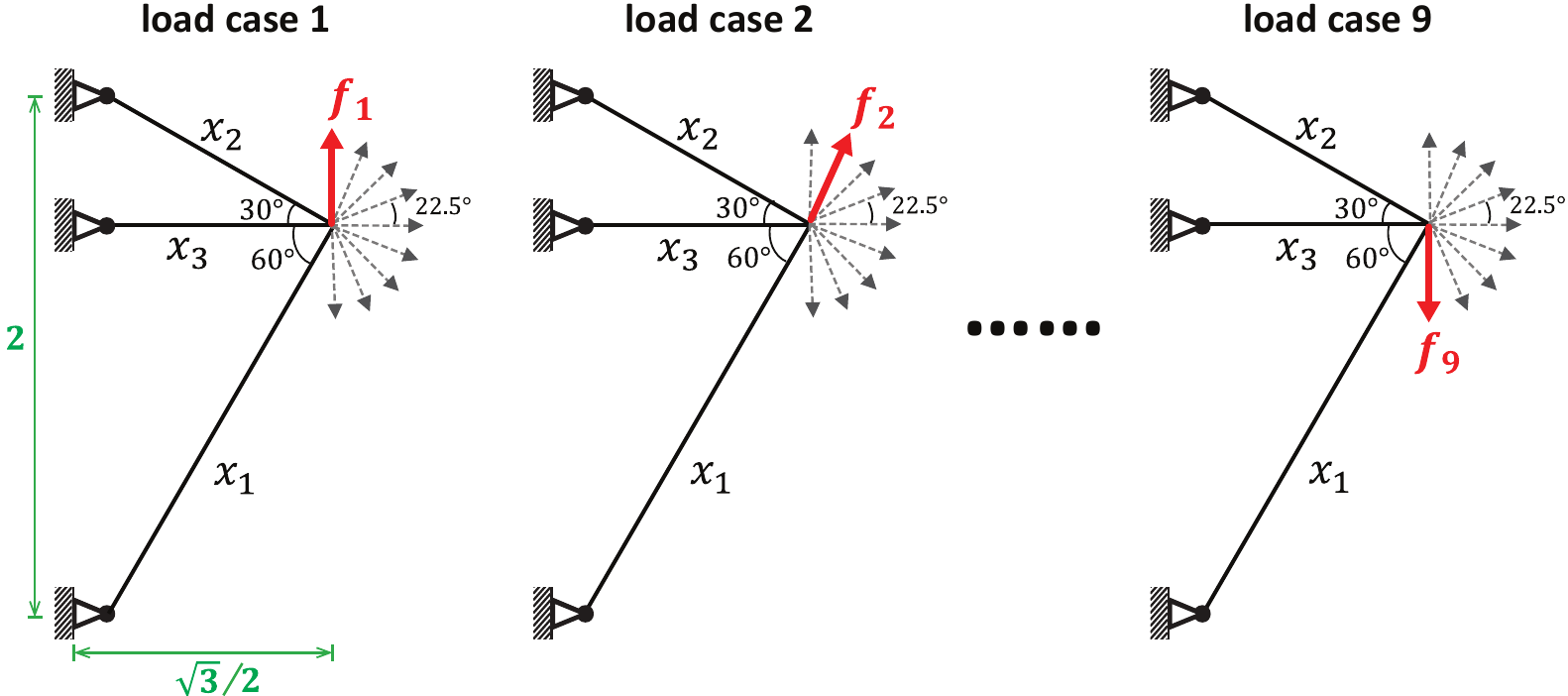}
	\caption{Three-bar truss structure under 9 load cases: initial ground structure, load and boundary conditions.}\label{Fig1}
\end{figure}
%
\begin{equation}\label{eq:demo_problem}
	\begin{aligned}
		&\min_{\boldsymbol{x}} C\left(\boldsymbol{x}\right)=\sum_{i=1}^{9}\textit{\textbf{f}}_i^{\,\,T}\textit{\textbf{u}}_i(\boldsymbol{x}) ,\\
		\text{s.t.}\quad & \sum_{e=1}^3 L^{(e)} x^{(e)}-V_{\text{max}}\leq0 \,,\\
		& 0<x^{\text{min}}\leq x^{(e)}\leq x^{\text{max}},\,e=1,...,3 ,\\
		& \text{with } \textit{\textbf{u}}_i(\boldsymbol{x})=\textit{\textbf{K}}\left(\boldsymbol{x}\right)^{-1}\textit{\textbf{f}}_i,\, i=1,...,9\,.\end{aligned}
\end{equation}
We set the volume constraint as $V_{\text{max}}=0.1$, and take initial guess as $x_0^{(e)}=V_{\text{max}}/\sum_{e} L^{(e)}=0.0278$,  for $e=1,2,3$. In addition, we choose the lower and upper bounds to be $x_{\text{min}}^{(e)}=10^{-8}x_{0}^{(e)}$  and $x_{\text{max}}^{(e)}=10^{4}x_{0}^{(e)}$, respectively. The initial move limit is taken as $\text{move}=10^{-1} x_0^{(e)}$. The tolerance for convergence of the optimization is $10^{-8}$. For the stochastic algorithm, we choose the sample size as $n=6$ (i.e., 6 stochastic loads), and select a new sample at each iteration.

\begin{figure}[!htbp]
	\centering
	\includegraphics[width=6.9in]{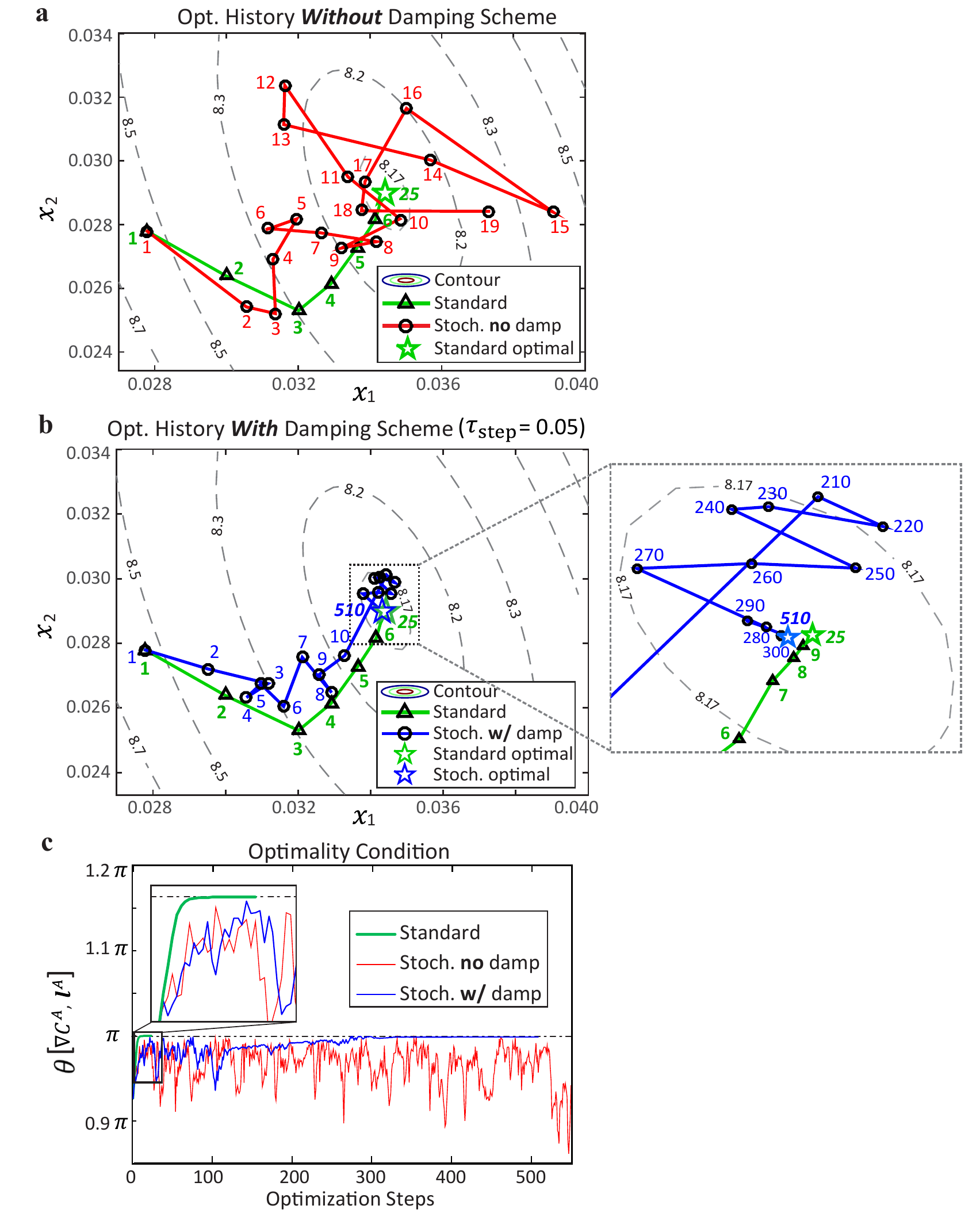}
	\caption{Illustration of the damping scheme in the three-bar truss example: contour plots of the objective function with the optimization history of $x^{(1)}$ and $x^{(2)}$ for (a) the standard algorithm and the randomized algorithm without a damping scheme; (b) the standard algorithm and the randomized algorithm with the proposed damping scheme. (c) The angle between the reduced gradient vectors of the objective function and the volume constraint.}\label{Fig2}
\end{figure}

We use the proposed stochastic algorithm with and without the damping scheme. The contour plots of the objective function with the optimization history of $x^{(1)}$ and $x^{(2)}$ 
($x^{(3)}$ can be computed from the volume constraint because, in practice, the sum of volumes will always be equal to $V_{\text{max}}$)
 for both cases are shown in Figs. \ref{Fig2}a and \ref{Fig2}b. The optimization history from the standard algorithm is plotted for reference. The standard algorithm obtains the global minimum of the given problem within 25 steps where $\boldsymbol{x}^*=\left[x^{(1)},x^{(2)},x^{(3)}\right]^T=\left[0.0344,0.0290,0.0132\right]^T$ and $C\left(\boldsymbol{x}^*\right)=8.1666$.
The stochastic algorithm without damping does not converge, and the updates become ineffective after roughly 10 steps (Fig. \ref{Fig2}a). 
The stochastic algorithm with our damping scheme converges
to $\widehat{\boldsymbol{x}}^S=\left[0.0343,0.0290,0.0134\right]^T$ with $C(\widehat{\boldsymbol{x}}^S)=8.1667$. The results are summarized in Table \ref{tbl6}. This is close to the optimal solution obtained by the standard algorithm. For this small problem, we use $n_\text{step}=10$, $Tol_\text{step}=0.05$, and $\gamma=2$. 
Since the solution of this example is relatively trivial, the stochastic algorithm with our damping scheme converges slower than the standard algorithm. For more complicated and realistic problems (e.g., examples in Section 4), the convergence rates for the stochastic and standard algorithms are comparable. 

\begin{table}[h!]
\centering
\captionsetup{justification=centering}
\caption{Solution of the 3-bar truss of Figure \ref{Fig1} using the standard GSM and \\the GSM with stochastic sampling and damping.}
\label{tbl6}
\begin{tabular}{c c c}
	\hline
	 & Standard GSM & Randomized\\
	&   & GSM w/ damping \\
	\hline
	 $x^{(1)}$ & 0.0344 & 0.0344 \\
	 $x^{(2)}$ & 0.0290 & 0.0290 \\
	 $x^{(3)}$ & 0.0132 & 0.0134 \\
 \hline 	
 $C^{*}$ & 8.1666 & 8.1667 \\
 \hline
\end{tabular}
\end{table}


To verify that the solution from the stochastic algorithm with damping converges to a KKT (Karush-Kuhn-Tucker) point (optimal solution in case of the GSM), we examine the angle $\theta^\mathcal{A}$ between $\nabla C^\mathcal{A}$ (reduced gradient vector of objective function) and $\boldsymbol{L}^\mathcal{A}$ (reduced gradient vectors of volume constraint) for the standard algorithm, and stochastic algorithm with and without damping, as shown in Fig. \ref{Fig2}(c). The solution is a KKT point if $\theta^\mathcal{A}=\pi$. For the standard algorithm, $\theta^\mathcal{A}$ converges quickly to $\pi$. For the stochastic algorithm with damping, $\theta^\mathcal{A}$ gradually converges to $\pi$ (the case without damping does not converge). Hence, we numerically show that the solution from the stochastic algorithm with damping converges to the optimal solution in the truss optimization framework.

\subsection{Overview of algorithmic parameters for randomized optimization}

This subsection summarizes the important parameters that are used in the randomized optimization framework, with some comments on the possible range of values to be used in practice.

\paragraph{Sample size}
In the proposed randomized optimization framework, 
the larger the number of sample load cases $n$, the more accurate the estimate of the compliance will be.
However, as the number of the  load cases $n$ increases, the associated computational cost also increases, because in each optimization step we need to solve a system of equations $n$ times. Thus, we need to balance the accuracy of the estimates and the computational complexity. Typically, the results from the stochastic algorithm are relatively insensitive to the number of sample load cases. Indeed, with a small number of sample load cases $(n \ll m)$ we obtain comparable (in terms of topology and compliance value) or sometimes even better (in terms of compliance value) solutions compared to those from the standard algorithm. The reason could be that the estimated gradient is already in the descent direction for the given sample size and the optimization steps are effective overall. To provide some insight into this discussion, a study is conducted using different sample sizes in Section 4.1.  
For the examples in the remainder of this paper, we choose the number of load cases $n=6$, unless otherwise stated, to demonstrate the accuracy and computational efficiency of the randomized optimization framework.
%

\paragraph{Frequency to select a new sample}
The frequency to select a new sample (or the number of optimization steps with fixed randomized sample) influences the convergence rate of the optimization. A frequency that is too low could result in convergence to a configuration that reflects a specific set of random loads, which are far away from the original load cases and the steps may be biased, ultimately leading to slower convergence. In this work, we select the random sample every step, i.e., $n_s=1$.

\paragraph{Filter parameters}
Applying the discrete filter in the GSM reduces the number of redundant bars and helps to limit the effects of the stochastic estimates for the randomized algorithm. In this work, we choose a relatively small filter size $\alpha_f=10^{-4}$ and effectively remove truss members when their cross-sectional areas have remained below $\alpha_f$ for $n_f=10$ cumulative steps. This provides additional consideration for the proposed stochastic sampling technique. Further details and variation of the filter parameters can be found in the paper by Ramos Jr. and Paulino \cite{Ramos2016}


\paragraph{Parameters in the damping scheme}
The damping scheme introduced in Section \ref{dampingscheme} has four parameters, the effective step ratio ($R$), the sample window size ($n_\text{step}$), the tolerance ($\tau_\text{step}$), and the scale factor ($\gamma$). The parameter choices in this damping scheme are crucial to the quality of the final solution. The size of the sample window (number of monitored steps, $n_\text{step}$) relates to the accuracy of the effective step ratio.
The effective step ratio becomes more accurate if we evaluate over a large window size. However, if $n_\text{step}$ is too large, it slows down convergence. In practice, we have found that it is sufficient to use $n_\text{step}=100$ for problems containing more than one thousand design variables. The number of steps after which we start to dampen the allowable moving limit is typically chosen to be the same as the window size.

The tolerance for the effective step ratio, $\tau_\text{step}$, serves as a threshold to determine when updates become ineffective. Therefore, the choice of $\tau_\text{step}$ affects the rate of convergence and sometimes the quality of the solutions. In general, a loose tolerance leads to faster convergence (reduce move limit more frequently) at the expense of the quality of design (in terms of the compliance value). One must balance the quality of the results and the convergence rate. The GSM is more sensitive to $\tau_\text{step}$ than the density-based method. This is demonstrated in Section \ref{section:2dexample}. Therefore, a stricter tolerance is needed for the GSM. In practice, we choose $\tau_\text{step}=0.05$ for the GSM and  $\tau_\text{step}=0.1$ for the density-based method.
For the moving limit scale factor $\gamma$, we have found that $\gamma=2$ is typically a good choice.


%% file: Section4-NumericalExamples_v1.tex
\section{Numerical examples}

We present several numerical examples in both two and three dimensions to demonstrate the effectiveness as well as the computational efficiency of the proposed randomized algorithm for topology optimization. Both density-based method and GSM are used. The first two examples (Section \ref{section:2dexample}) investigate the sensitivity of the density-based method and the GSM to the tolerance $\tau_\text{step}$ (see Eq. \eqref{eq:stepev}). Example 2 (Section \ref{section:2dDensity}) further shows the relation between sample size and quality of the optimized design. Example 3 (Section \ref{section:2dGSM}) illustrates the effect of the discrete filter in the GSM for the proposed randomized algorithm. The last two examples (Sections \ref{section:3dDensity} and \ref{section:3dGSM}) in 3D demonstrate the capability of the proposed algorithm to create practical structural designs at greatly reduced computational cost.

To quantify the computational cost of the standard and stochastic optimization algorithms, we define $N_\text{solve}=n\times N_\text{step}$ as the total number of linear systems of equations to solve in the optimization process, where $N_\text{step}$ is the number of optimization steps. This is a measure of the computational efficiency of an optimization formulation. The optimization process is considered converged if the current step size (bounded by the move limit) is below a prescribed tolerance $\tau_\text{opt}$ for the optimization process, that is, $||\boldsymbol{x}_k-\boldsymbol{x}_{k-1}||<\tau_\text{opt}$.

For the density-based method (continuum), we incorporate the proposed technique into the computer program PolyTop \cite{Talischi2012a} in 2D and the topology optimization code in 3D \cite{liu2014}. For plotting 3D continuum results, we utilize TOPSlicer \cite{zegard2016}. All the problems are initialized as follows. The initial guess is taken as $\rho_0^{(e)}=V_{\text{max}}/M$, where $M$ is the number of elements in the finite element mesh. The convergence tolerance is $\tau_\text{opt}=10^{-2}$; the initial move limit is chosen as $\text{move}=0.05$; the damping factor for the OC update scheme is $\eta=0.5$. The penalization factor starts at $p=1$ and gradually increases to $p=3$ (for 2D problems) or $4$ (for 3D problems).

For the GSM (discrete), we generate initial ground structures (without overlapped bars) using the collision zone technique from references \cite{zegard2014,zegard2015a} and plot final topologies in 3D using the program GRAND3 \cite{zegard2015a}. The initial guess of the design variables is taken as $x_0^{(e)}=V_{\text{max}}/\sum_{e=1}^{M}L^{(e)}$; the convergence tolerance is $\tau_\text{opt}=10^{-8}$; the initial move limit is chosen as $\text{move}=x_0^{(e)}\times10^{4}$; the damping factor for the OC update scheme is $\eta=0.5$. For the case that the discrete filter is used in the GSM, we use $n_f=10$ and $\alpha_f=10^{-4}$ during the optimization process unless otherwise stated; the lower and upper bounds on the design variables are $x_{\text{min}}=0$ and $x_{\text{max}}=10^4x_0$ (unbounded in practical terms).   For the standard GSM (without the discrete filter), we apply a cut-off value $10^{-2}$ that defines the final structure at the end of the optimization \cite{christensen2009}. The lower and upper bounds are defined by $x_{\text{min}}=10^{-2}x_0$  and $x_{\text{max}}=10^{4}x_0$, respectively. For all results in the GSM, we remove unstable nodes and floating bars and then check the final topologies to ensure that they are at global equilibrium--detailed explanation can be found in references \cite{Ramos2016,zhang2016a}.

In the above two optimization methods (continuum and discrete types), the stochastic algorithm uses the following parameters. Unless otherwise stated, the sample size is chosen to be $n=6$; in the damping scheme, the window size is taken to be $n_\text{step}=100$, we assign the step size reduction factor as $\gamma=2$ and the tolerance for the effective step ratio as $\tau_\text{step}=0.1$ for the density-based method and $\tau_\text{step}=0.05$ for the GSM. Let $\boldsymbol{\rho}^{*}$ and $\boldsymbol{x}^{*}$ represent the optimal solutions of the standard formulations in \eqref{eq:density_obj} and \eqref{eq:ground_obj}, $\widehat{\boldsymbol{\rho}}^{S}$ and $\widehat{\boldsymbol{x}}^{S}$ represent the optimal solutions obtained from the stochastic algorithm for density-based and ground structure methods. To evaluate the quality of the solutions, in the case of the stochastic algorithm, we present the true values of the objective function $C(\widehat{\boldsymbol{\rho}}^{S})$ and $C(\widehat{\boldsymbol{x}}^{S})$ at the approximated solutions $\widehat{\boldsymbol{\rho}}^{S}$ and $\widehat{\boldsymbol{x}}^{S}$ (instead of their estimators $\widehat{C}^S(\widehat{\boldsymbol{\rho}}^{S})$ and $\widehat{C}^S(\widehat{\boldsymbol{x}}^{S})$) and compare them with those obtained from the standard algorithm $C({\boldsymbol{\rho}}^{*})$ and $C({\boldsymbol{x}}^{*})$. The relative difference is defined as $\Delta C= \left(C(\widehat{\boldsymbol{x}}^{S})-C({\boldsymbol{x}}^{*})\right)/C({\boldsymbol{x}}^{*})$.



\subsection{Two-dimensional box domain with 108 load cases}\label{section:2dexample}
\begin{figure}[!htbp]
	\centering
	\includegraphics[width=3.5in]{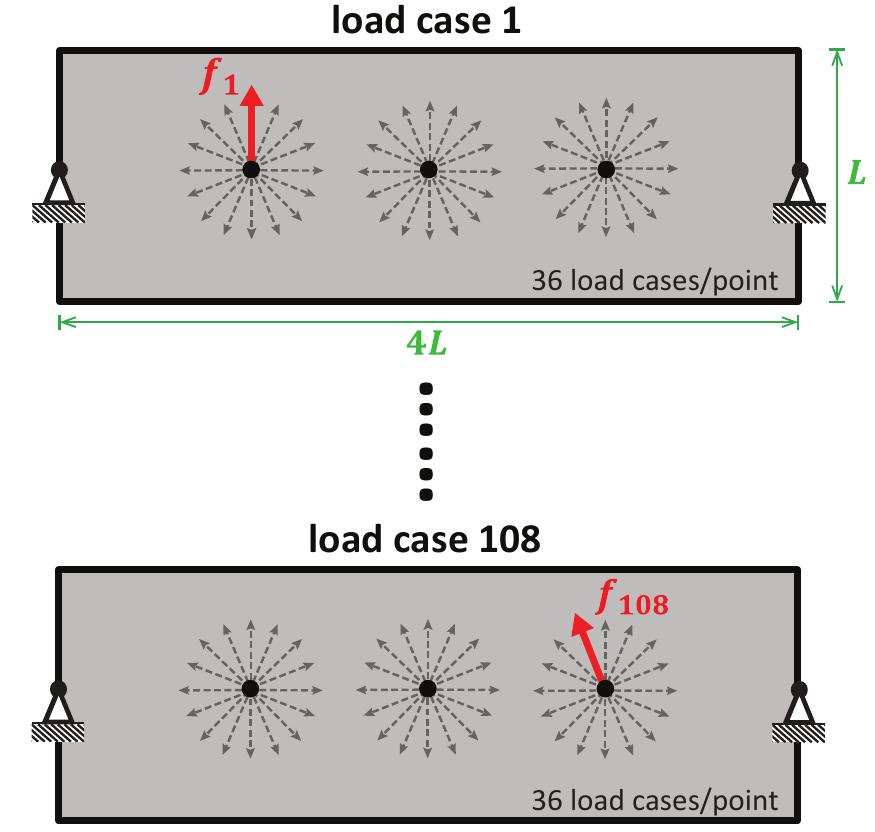}
	\caption{Two-dimensional box domain with load and support conditions. A total of 108 equal-weighted load cases are applied at three given points with each point having 36 load cases applied from $0^{\circ}$ to $350^{\circ}$ (dotted arrows are the schematic illustrations of non-active load cases).}\label{Fig3}
\end{figure}

We present a two-dimensional (2D) topology optimization problem whose design domain and boundary conditions are shown in Fig. \ref{Fig3}. A total of 108 equal-weighted load cases are applied at three given points, with each point having 36 load cases applied along different angles (from $0^{\circ}$ to $350^{\circ}$). In this section, both the density-based and the ground structure methods are used.
\subsubsection{Example 1: Continuum topology optimization with density-based method}\label{section:2dDensity}

\begin{figure}[!htbp]
	\centering
	\includegraphics[width=6.1in]{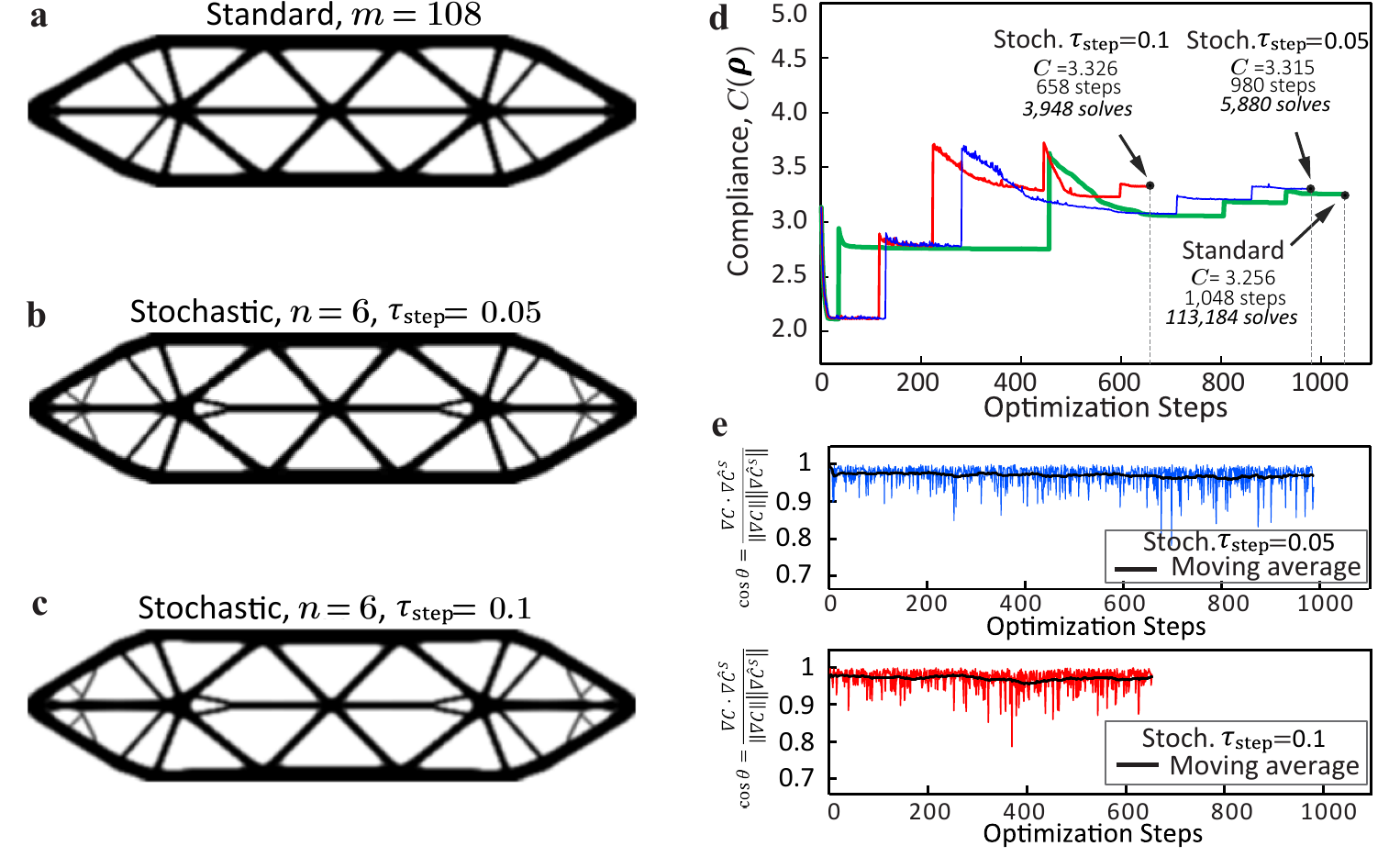}
	\caption{Results for Example 1 using the density-based method with 25,600 quadrilateral elements and 52,002 degrees of freedom (DOFs). (a) The optimized topology obtained by the standard algorithm \cite{bendsoe2003}; (b) the optimized topology obtained by the stochastic algorithm with $n=6$ and $\tau_\mathrm{step}=0.05$ (one representative trial); (c) the optimized topology obtained by the stochastic algorithm with $n=6$ and $\tau_\mathrm{step}=0.1$ (one representative trial); (d) the convergence of the compliance for above cases; (e) cosine of the angle between the gradient $\nabla C_x$ and the estimated gradient $\nabla \widehat{C}_x^{S}$ for the stochastic cases ($\tau_\mathrm{step}=0.05$ and $\tau_\mathrm{step}=0.1$) to verify if their directions are aligned ($\cos\theta\approx 1$). }\label{Fig4}
\end{figure}

Using the density-based method, we demonstrate the reduction of the computational cost by means of the stochastic algorithm. We further investigate the sensitivity of the final optimized topologies to the tolerance $\tau_\text{step}$ in the damping scheme and the sample size $n$. Since the final topology from the standard algorithm is symmetric both horizontally and vertically, we enforce symmetry of the topologies in the stochastic case by enforcing horizontal and vertical symmetry of the distribution of the density field \cite{Talischi2012a}. A total number of 25,600 quadrilateral elements are used to discretize the domain which gives 52,002 degrees of freedom (DOFs). The linear density filter that defines the solid-void boundary takes the radius of $1$ (see Section 1.1). For comparison purposes, the final topology obtained from the standard topology optimization is shown in Fig. \ref{Fig4}(a). The final topology has $C(\boldsymbol{\rho}^{*})=3.257$ and converges after 1048 steps. In each optimization step, we solve 108 linear systems (corresponding to 108 load cases), which leads to a total $113,184$ solves. Since the continuation method is used for the penalization factor, $p$, the jumps in the compliance correspond to $p=1.5, 2.0, 2.5, 3.0$ (initially $p=1.0$).

To investigate the sensitivity to $\tau_\text{step}$, we consider $\tau_\text{step}=0.1$ and $\tau_\text{step}=0.05$,
and the results are compared with that from the standard algorithm \cite{bendsoe2003}. For both cases, the number of sample load cases used is $n=6$. Figures \ref{Fig4}(b) and (c) show the optimized topologies for the stochastic algorithm for a single representative trial (one trial is one run of the numerical experiment) for each $\tau_\text{step}$, and Fig. \ref{Fig4}(d) shows the convergence histories of the objective function for all cases for the representative trial. Since the sample load cases are generated randomly at each step, the final optimized topology and its compliance vary with each trial. Therefore, the associated results in Table \ref{tbl1} are averaged over 5 trials.  For this example, the standard algorithm leads to a lower compliance and simpler topology. However, the stochastic algorithm for both tolerances uses substantially fewer solves for final topologies similar to the one obtained with the standard algorithm. For $\tau_\text{step}=0.1$, the number of linear systems to solve is 27 times fewer than for the standard algorithm (113,184 solves vs. 4170 solves), and the convergence of the optimization is more rapid.
The final topologies obtained for the two tolerances and the compliance for each case (see Fig. \ref{Fig4})
suggest that the tolerance in the damping scheme has a minor influence on the final results in the density-based method. The relative differences between the compliance for the standard algorithm and those for the stochastic algorithm are $1.79\%$ for $\tau_\text{step}=0.05$ and $2.45\%$ for $\tau_\text{step}=0.1$. It seems that smaller $\tau_\text{step}$ leads to slower convergence but a slightly better compliance: the optimization with $\tau_\text{step}=0.05$ takes an average 1052 steps to converge while the one with $\tau_\text{step}=0.1$ takes an average 695 steps. Therefore, the latter is more computationally efficient with only $N_\text{solve}=4170$ compared with $N_\text{solve}=6312$ for the former one.

Next, we check the quality of the estimated gradients by showing that the estimated negative gradient of the objective function is a descent direction. As shown in Fig. \ref{Fig4}(e), we plot the cosine of the angle ($\theta$) between the gradient and the estimated gradient for each of the two tolerances for one representative trial. The moving averages (over 50 steps) of $\cos\theta$, in both cases, are close to 1.0.

\begin{figure}[!htbp]
	\centering
	\includegraphics[width=6.5in]{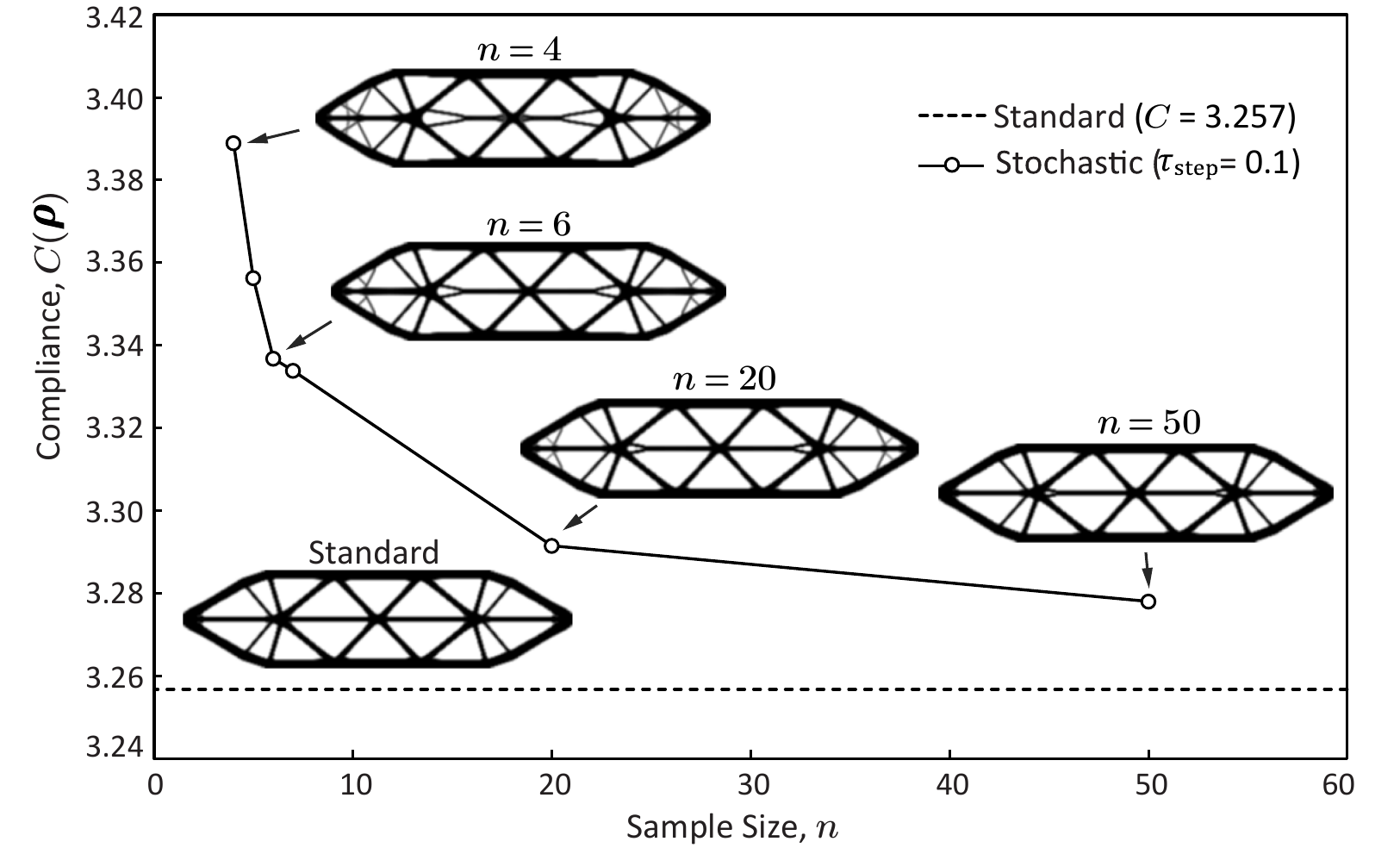}
	\caption{Study of sample sizes ($n=4,5,6,7,20,50$) versus the resulting final compliance (or end-compliance) using the stochastic algorithm in Example 1. The final topologies (from representative trials) are included.}\label{Fig5}
\end{figure}

The next study demonstrates the influence of sample sizes on optimization results and the reduction of the computational cost by using the stochastic algorithm. In this study, $\tau_\text{step}=0.1$. Figure \ref{Fig5} gives the compliances for sample sizes $n=4,5,6,7, 20, 50$ and final topologies (from one representative trial). Each data point in Fig. \ref{Fig5} is obtained by averaging 5 trials, and the data are summarized in Table \ref{tbl1}. Several observations can be made based on Fig. \ref{Fig5} and Table \ref{tbl1}. For one, the stochastic algorithm leads to similar optimal topologies and compliances compared with the standard algorithm. For $n=4$, $N_\text{solve}$ is significantly smaller than for the standard algorithm by a factor of 45 on average. Moreover, the compliance improves as we increase $n$, indicating that larger $n$ offers better estimation during the optimization and ultimately yields stiffer optimal structures. However, $N_\text{solve}$ also increases as we use more sample load cases. From the optimized structures and compliances, it seems that $n=6$ is sufficient for this problem with greatly reduced computational cost. The compliance differs from the standard algorithm by only $2.45\%$.
Table \ref{tbl1} shows that the average cosine of angles, $\cos\bar{\theta}$, between the gradient and estimated gradient ranges from $0.959$ and $0.996$ for various sample sizes.
Using the stochastic algorithm, we can almost fully recover the original optimal results by either increasing the sample size (e.g., $n=50$) or choosing smaller $\tau_\text{step}$, as both methods lead to highly accurate designs.


%
\begin{table}[h!]
	\centering
	\captionsetup{justification=centering}
	\caption{Results for Example 1 (density-based), averaged over 5 trials.}
	\label{tbl1}
	\begin{tabular}{c c c c c c c c c c }
		\hline
		Density-based & \parbox{13mm}{$C(\boldsymbol{\rho}^{*})$} &  $C(\widehat{\boldsymbol{\rho}}^{S})$ & $\Delta C$ & $\tau_\mathrm{step}$ & $\cos\bar{\theta}$& $n$  & $N_\text{step}$ & $N_\text{solve}$ \\
		\hline
		 Standard        & 3.257 & - & -& - & - & - & 1,048 & 113,184  \\
		Stoch. $\tau_\mathrm{step}=0.05$ & - & 3.315 & 1.79\%&  $0.05$  & $0.971$& 6 & 1,052 & 6,312\\
		\textbf{Stoch. $\bm{\tau_\mathrm{step}=0.1}$ }& - & \textbf{3.337} & \textbf{2.45\%}& $\bm{0.1}$  & $\bm{0.971}$& \textbf{6} &   \textbf{695} & \textbf{4,170} \\
		Stoch. $n=4$ & - & 3.389 & 4.05\%&  $0.1$  & $0.959$& 4 & 618 & 2,472 \\
		Stoch. $n=5$ & - & 3.356 & 3.05\%& $0.1$ & $0.967$& 5 & 684 & 3,420 \\
		Stoch. $n=7$ & - & 3.334 & 2.36\%& $0.1$  & $0.976$& 7 & 684 & 4,788 \\
		Stoch. $n=20$ & - & 3.291 & 1.05\%&  $0.1$  & $0.991$& 20 & 838 & 16,760 \\
		Stoch. $n=50$ & - & 3.278 & 0.65\%& $0.1$  & $0.996$& 50 & 850 & 42,500 \\\hline 	
	\end{tabular}
\end{table}

\subsubsection{Example 2: Truss topology optimization with ground structure method}\label{section:2dGSM}

\begin{figure}[!htbp]
	\centering
	\includegraphics[width=6.5in]{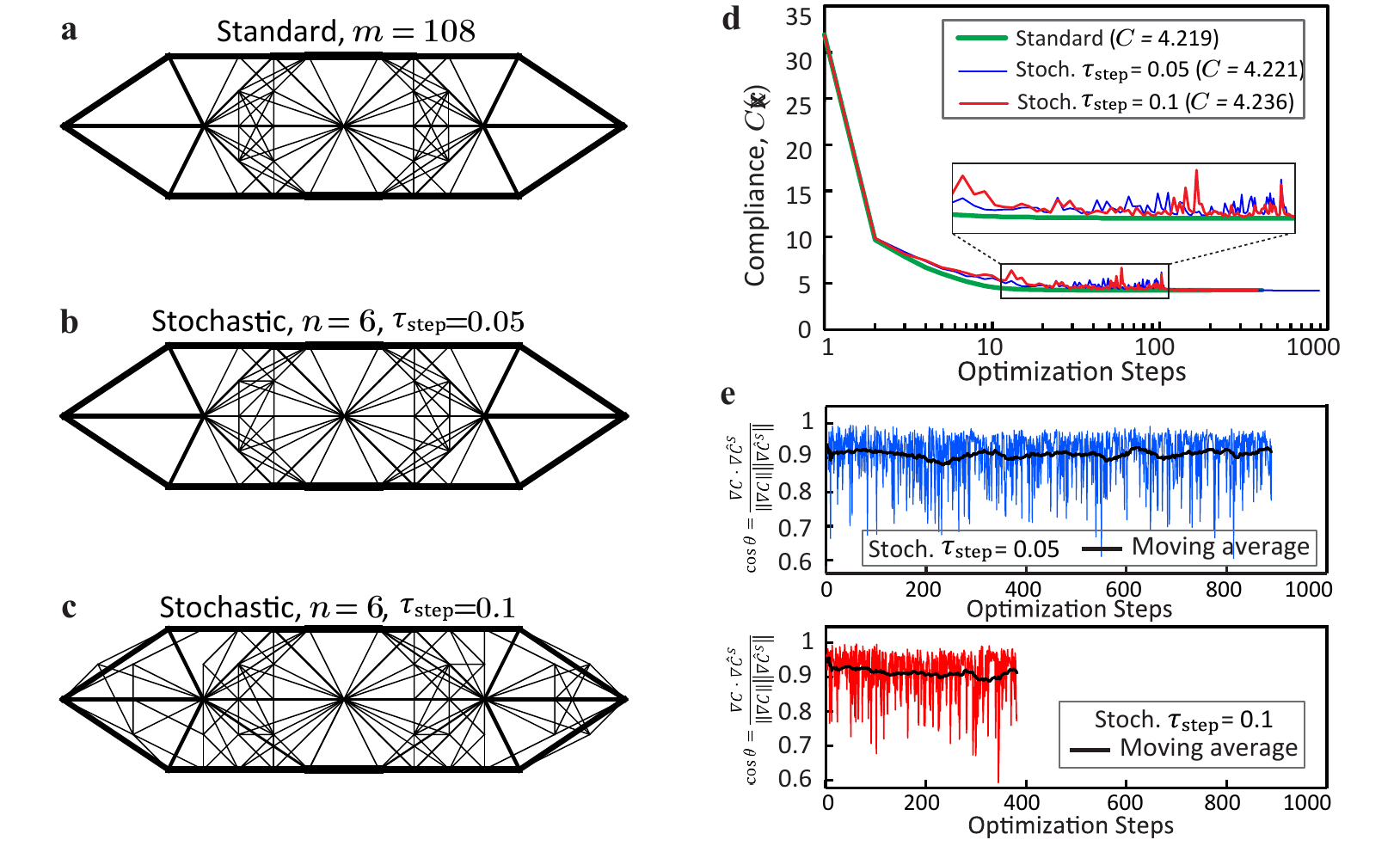}
	\caption{Results for the GSM (without the discrete filter) for Example 2 with a full-level ground structure ($16\times 4$ grid) and 2,196 non-overlapped bars. (a) The optimized structure obtained by the standard algorithm; (b) the optimized structure obtained by the stochastic algorithm with $n=6$ and $\tau_\mathrm{step}=0.05$; (c) the optimized structure obtained by the stochastic algorithm with $n=6$ and $\tau_\mathrm{step}=0.1$; (d) the convergence of the compliance for all above cases; (e) cosine of $\theta$ between the gradient $\nabla C_x$ and the estimated gradient $\nabla \widehat{C}_x^{S}$ for the stochastic algorithm ($\tau_\mathrm{step}=0.05$ and $\tau_\mathrm{step}=0.1$). }\label{Fig6}
\end{figure}

As example 2, we study again the problem presented in Fig. 3, but
this time using the ground structure method. We demonstrate that our approach greatly reduces the total number of linear solves. In addition, we investigate the sensitivity of the results to $\tau_\mathrm{step}$ as well as the influence of the discrete filter on final solutions of the stochastic algorithm.  We use a full-level ground structure ($16\times 4$ grid) with 2,196 non-overlapped bars to discretize the domain \cite{zegard2014,zegard2015a}. The optimal topology and the convergence of the compliance for the standard algorithm are shown in Figs. \ref{Fig6}(a) and (d). The optimal compliance for the standard algorithm, $C\left(\boldsymbol{x}^*\right)=4.219$, is obtained in 406 steps and $N_\text{solve}=43308$. For the stochastic algorithm, we compare two cases using different tolerances in the damping scheme, $\tau_\mathrm{step}=0.05$ and $\tau_\mathrm{step}=0.1$. The results are summarized in Table \ref{tbl2}, averaged over 5 trials. The optimized structures and the convergence of the compliance for a single representative trial are shown in Figs. \ref{Fig6}(b)-(d).

In contrast to the results for the density-based method in Section 4.1.1, the results for the GSM indicate that the choice of $\tau_\mathrm{step}$ has a significant impact on the optimized structure. Although $\tau_\mathrm{step}=0.05$ results in a larger number iterations/optimization steps and a larger number of linear system solves, $N_\text{solve}$, compared with $\tau_\mathrm{step}=0.1$, its final structure is simpler, has slightly lower compliance, and is similar to the structure obtained by the standard algorithm. This suggests that $\tau_\mathrm{step}=0.05$ is a better choice for this example. The convergence rate for the stochastic algorithm is about the same as for the standard algorithm. Since the standard optimization problem for the GSM is convex, its solution is the global minimum; hence, the compliances obtained with the stochastic algorithms will be larger than or equal to those obtained with the standard algorithm. However, the relative differences in compliance are very small, only $0.09\%$ and $0.35\%$ in average, and the stochastic algorithm achieves these results with considerably less computational effort ($N_\text{solve}=5627$ and $2,322$ on average for the stochastic algorithm versus $N_\text{solve}=43,848$ for the standard algorithm). Symmetry was not enforced for the GSM in the present study.
To show the accuracy of estimated gradient of compliance for both tolerances, we plot $\cos\theta$ between the gradient and the estimated gradient for one representative trial in Fig. \ref{Fig6}(e). The moving averages (over 50 steps) of $\cos\theta$ typically stay above 0.9, indicating that the estimated gradients of both stochastic cases are quite accurate and lead to effective optimization steps.

\begin{figure}[!htbp]
	\centering
	\includegraphics[width=6.5in]{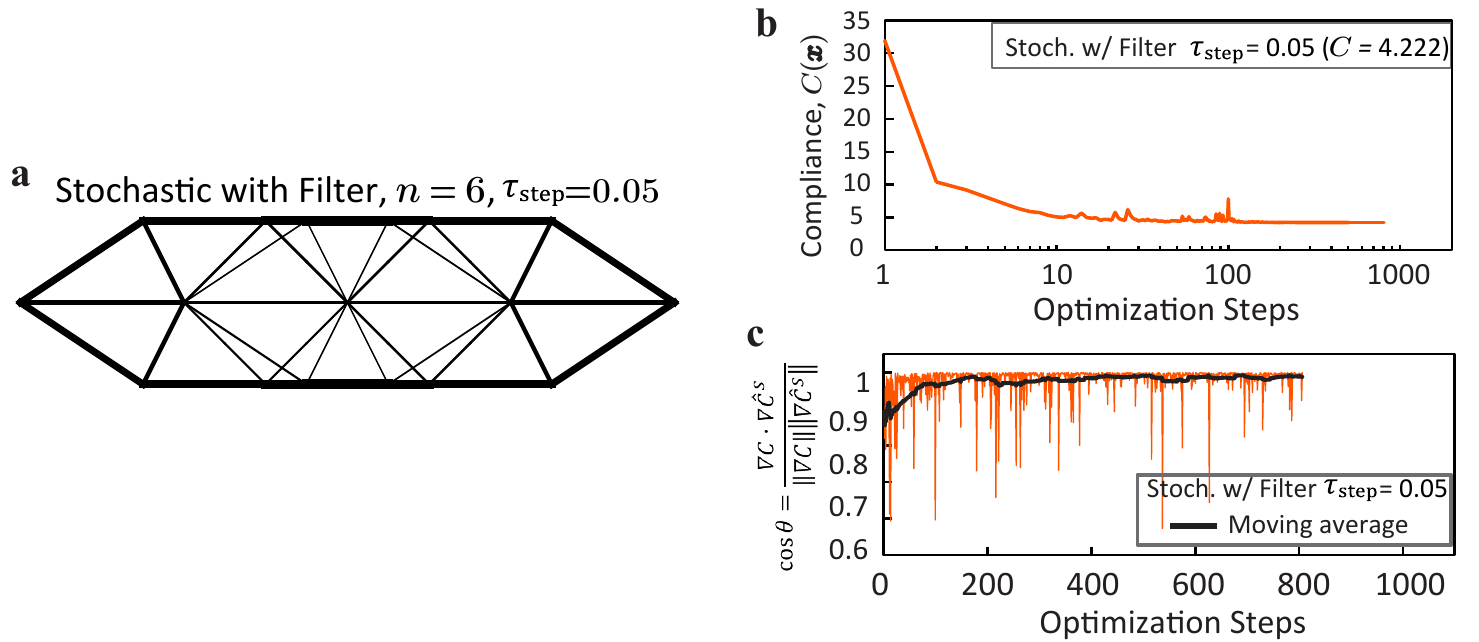}
	\caption{Results for the GSM with the discrete filter for Example 2. (a) The optimized structure obtained from stochastic algorithm with $n=6$, $\tau_\mathrm{step}=0.05$, and $\alpha_{f}=0.0001$; (b) the convergence of the compliance; (c) cosine of $\theta$ between the gradient $\nabla C_x$ and the estimated gradient $\nabla \widehat{C}_x^{S}$ for the stochastic case. }\label{Fig7}
\end{figure}

Next, we also include the discrete filter \cite{Ramos2016} in the stochastic algorithm and study its influence on the optimized results. Based on previous observations, we choose $\tau_\mathrm{step}=0.05$ and $n=6$. We use the following parameters for the discrete filter (see Section 2.3):  $n_f=10$ and $\alpha_f=0.0001$. Figure \ref{Fig7}(a) shows the final topology obtained by the stochastic algorithm with the discrete filter in one representative trial (cf. Fig. \ref{Fig6}(b), which was obtained without the filter). The convergence of the compliance and $\cos\theta$ in each step are shown in Figs. \ref{Fig7} (b)-(c). The results, averaged over five trials, are summarized in Table 2. {The discrete filter combined with the stochastic algorithm leads to a simpler final topology. From Fig. \ref{Fig7}(c), the direction of the estimated gradient seems to be more accurate than the ones without the discrete filter (Fig. \ref{Fig6}(e)), which is confirmed by $\cos\bar{\theta}=0.978$. This indicates that the removal of some non-useful members helps to limit the collateral effects of the stochastic estimates. As compared to the standard algorithm, the discrete filter combined with the stochastic algorithm not only leads to much less number of linear system solves ($N_\text{solve}=5,442$ versus $N_\text{solve}=43,848$), but the size of the linear systems also keeps decreasing by the discrete filter, which further improves the computational efficiency. }

\begin{table}[h!]
	\centering
	\captionsetup{justification=centering}
	\caption{Results for Example 2 (GSM), averaged over 5 trials.}
	\label{tbl2}
	\begin{tabular}{c c c c c c c c c c c }
		\hline
		GSM & \parbox{13mm}{$C(\boldsymbol{x}^{*})$} &  $C(\widehat{\boldsymbol{x}}^{S})$ & $\Delta C$ & $\tau_\mathrm{step}$ & $\cos\bar{\theta}$& $n$  &  $N_\text{step}$ & $N_\text{solve}$ \\
		\hline
		Standard        & 4.219 & - & -& -& - & - & 406 & 43,848  \\
		Stoch. $\tau_\mathrm{step}=0.05$ & - & 4.222 & 0.09\%&  $0.05$  & $0.911
		$& 6& 938 & 5,627\\
		Stoch. $\tau_\mathrm{step}=0.1$ & - & 4.233 & 0.35\%& $0.1$  &  $0.912$& 6& 387 & 2,322 \\
		\textbf{Stoch. w/ Filter} $\bm{\tau_\mathrm{step}=0.05}$ & - & \textbf{4.222} & \textbf{0.09\%}& $\bm{0.05}$  &  $\bm{0.978}$& \textbf{6}& \textbf{907} & \textbf{5,442} \\  \hline 	
	\end{tabular}
\end{table}



\subsection{Example 3: Three-dimensional bridge design using the density-based method}\label{section:3dDensity}
\begin{figure}[!htbp]
	\centering
	\includegraphics[width=6.5in]{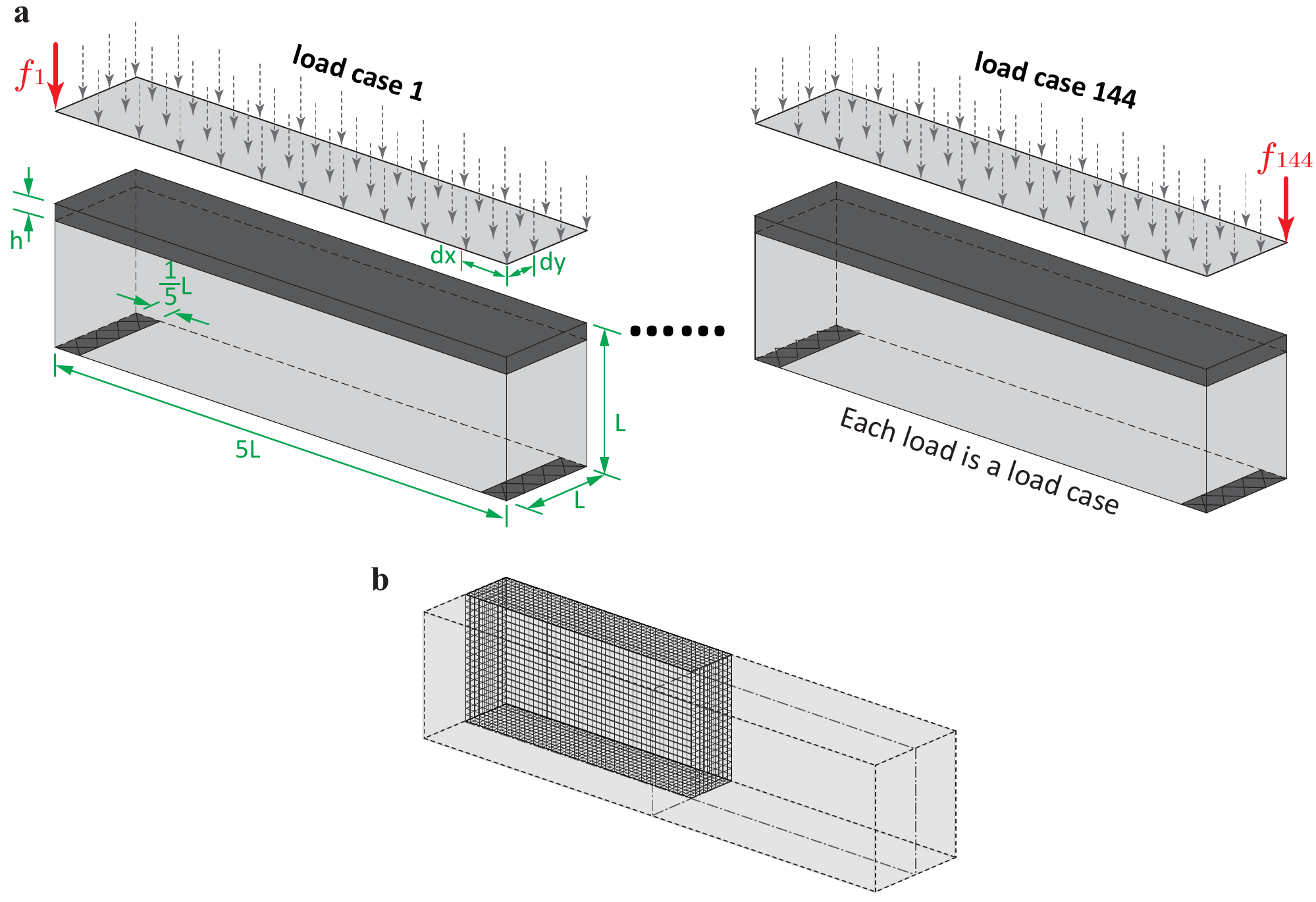}
	\caption{Three-dimensional bridge design with (a) the geometry, load and boundary conditions; (b) the quarter domain is modeled by a mesh with 10,000 brick elements and 35,343 DOFs.}\label{Fig8}
\end{figure}

\begin{figure}[!htbp]
	\centering
	\includegraphics[width=6.2in]{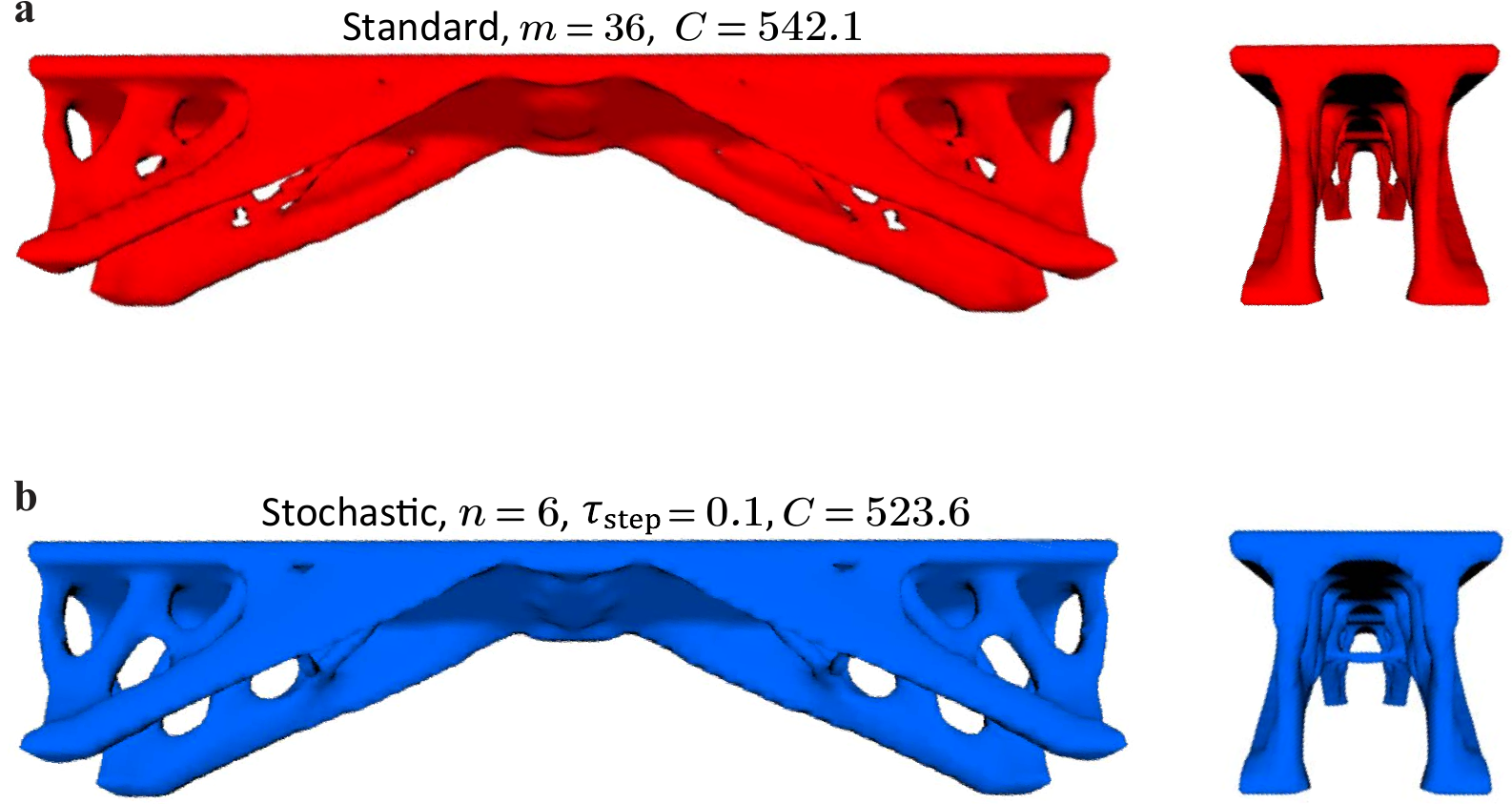}
	\caption{Optimized structures of the 3D bridge design obtained from (a) the standard algorithm; (b) the stochastic algorithm with $n=6$ and $\tau_\mathrm{step}=0.1$. }\label{Fig9}
\end{figure}

In this section, we demonstrate the quality of the design and the great reduction in computational work of the proposed randomized algorithm using a 3D bridge design in the non-convex, continuum topology optimization framework. The design domain, load, and boundary conditions are shown in Fig. \ref{Fig8}(a). A total of 144 equal-weighted load cases are applied at the bridge deck (non-designable layer). Based on the structural symmetry, we optimize a quarter of the domain, as shown in Fig. \ref{Fig8}(b), which reduces the number of load cases to $m=36$. We use 10,000 brick elements to discretize the quarter domain, resulting in 35,343 degrees of freedom (DOFs). We take $V_{\text{max}}=0.1\times M$ (where $M=10,000$) and use a quadratic density filter which takes the radius of $2.5$ (see Section 1.1). The penalization factor for the continuation scheme takes values $p=1, 2, 3, 4$. The stochastic case uses the following parameter values: the sample size is chosen to be $n=6$; the window size $n_\text{step}=100$, $\gamma=2$ and $\tau_\mathrm{step}=0.1$ (the effective step ratio is calculated in terms of $\boldsymbol{\rho}$).  Figure \ref{Fig9} shows the optimized structures obtained using the standard algorithm and the stochastic algorithm. The results are summarized in Table \ref{tbl3}.

The standard algorithm leads to final topology with $C(\boldsymbol{\rho}^{*})=542.1$ and converges in 968 steps. In each optimization step, we solve 36 linear systems (load cases), which leads to $N_\text{solve}=34,848$. Our stochastic algorithm, while offering a nearly identical topology as the standard algorithm, drastically reduces the computational cost from $34,848$ solves to $4,662$ solves and leads to even better compliance $C(\widehat{\boldsymbol{\rho}}^{S})=523.6$.
\begin{table}[h!]
	\centering
	\captionsetup{justification=centering}
	\caption{Results for Example 3 (density-based bridge design).}
	\label{tbl3}
	\begin{tabular}{c c c c c c c c c c c }
		\hline
		Study 1 & \parbox{13mm}{$C(\boldsymbol{\rho}^{*})$} &  $C(\widehat{\boldsymbol{\rho}}^{S})$ & $\Delta C$ & $\tau_\mathrm{step}$ &  $n$  &  $N_\text{step}$ & $N_\text{solve}$ \\
		\hline
		 Standard        & 542.1 & - & -& -& - &  968 & 34,848  \\
		Stochastic & - & 523.6 & -3.4\%& $0.1$  &  6& 777 & 4,662 \\  \hline 	
	\end{tabular}
\end{table}
\subsection{Example 4: Three-dimensional high-rise building design using the ground structure method}\label{section:3dGSM}

\begin{figure}[!htbp]
	\centering
	\includegraphics[width=6.5in]{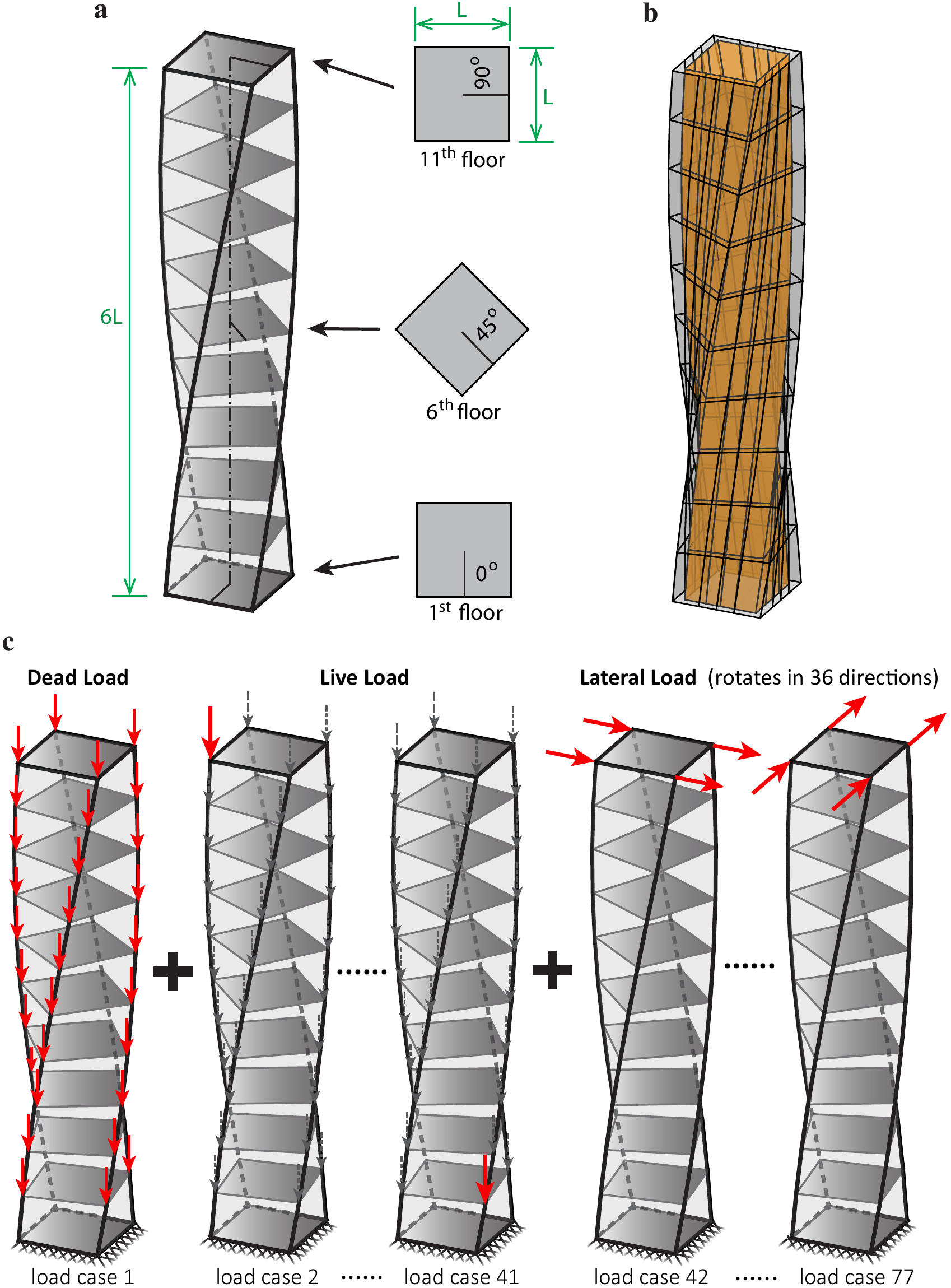}
	\caption{Twisting tower (inspired by the Canyan tower \cite{canyantower}): (a) geometry; (b) base mesh with a void zone in the middle; (c) load and boundary conditions. One dead load case, 40 live load cases, and 36 lateral load cases (the lateral load is applied at 4 corners on the top floor and rotating in 36 directions).}\label{Fig10}
\end{figure}

\begin{figure}[!htbp]
	\centering
	\includegraphics[width=6.5in]{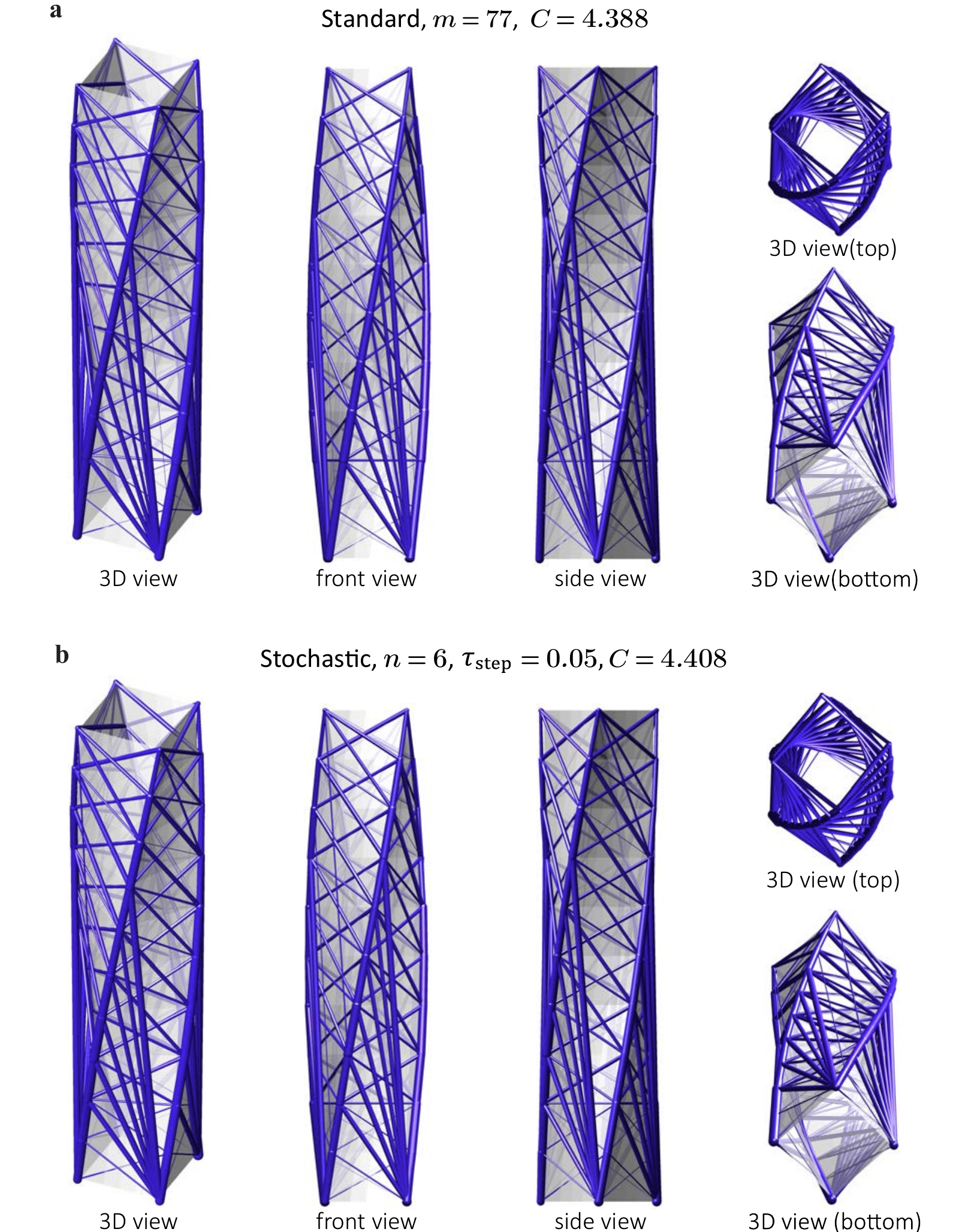}
	\caption{Optimized structures of the 3D twisting tower obtained from: (a) the standard algorithm; (b) the proposed stochastic algorithm with $n=6$ and $\tau_\mathrm{step}=0.05$.}\label{Fig11}
\end{figure}

To illustrate the effectiveness of the stochastic algorithm combined with the discrete filter on a practical engineering design, we optimize a simplified 3D high-rise building (twisting tower) in the truss topology optimization framework. The domain of the twisting tower, given in Fig. \ref{Fig10}, has 11 floors with the 1st floor fixed to the ground. The tower twists a full 90 degrees from its base to its crown. To obtain constructible structures, we use a $4\times4\times11$ grid to discretize the domain followed by the generation of a level 7 initial ground structure, containing 3,556 non-overlapping members \cite{zegard2015a}. The base mesh and the void zone are shown in Fig. \ref{Fig10}(b). As shown Fig. \ref{Fig10}(c), 77 equal-weighted load cases are applied to the building. Based on our previous studies, we choose $\tau_\mathrm{step}=0.05$ and $n=6$. We apply the discrete filter for both the standard and the stochastic algorithms; specifically, we choose $n_f=10$, a small filter value ($\alpha_f = 0.0001$) used during optimization, and a larger filter value ($\alpha_f = 0.001$) in the final step to control the resolution of the final topology \cite{Ramos2016,zhang2016}. Figure \ref{Fig11} shows the optimized structures obtained using the standard and stochastic algorithms. A summary of the results is provided in Table \ref{tbl4}.

The optimal compliance for the algorithm, $C\left(\boldsymbol{x}^*\right)=4.388$, is obtained in 382 optimization steps and a total of $N_\text{solve}=29,414$ linear solves. With the stochastic algorithm, we obtain a similar final structure as well as a compliance that is only 0.46\% higher than that obtained with the standard algorithm. These results are obtained at a drastically reduced computational cost, i.e., $N_\text{solve}=4,410$ linear solves. In addition to solving a less number of linear systems by using the stochastic algorithm, the discrete filter also reduces the size of the linear systems as the optimization proceeds because the use of the filter scheme reduces the number of design variables, the size of the stiffness matrices, and the size of the sensitivity vectors, which significantly decrease the CPU time and memory usage \cite{zhang2016}, contributing to great computational efficiency.


\begin{table}[h!]
\centering
\captionsetup{justification=centering}
\caption{Results for Example 4 (Tower design).}
\label{tbl4}
\begin{tabular}{c c c c c c c c c c c }
	\hline
	3D GSM & \parbox{13mm}{$C(\boldsymbol{x}^{*})$} &  $C(\widehat{\boldsymbol{x}}^{S})$ & $\Delta C$ & $\tau_\mathrm{step}$ & $\cos\bar{\theta}$& $n$  &  $N_\text{step}$& $N_\text{solve}$ \\
	\hline
	 Standard        & 4.388 & - & -& -& - & - & 382 & 29,414
	   \\
	Stochastic & - & 4.408 & 0.46\%&  $0.05$  & 0.937& 6& 735 & 4,410\\
 \hline 	
\end{tabular}
\end{table}

%% file: Section5-ConcludingRemarks_v1.tex
\section{Concluding remarks and future work}
This paper proposes an efficient randomized optimization approach for topology optimization that drastically reduces the enormous computational cost of optimizing practical structural designs under many load cases and produces high-quality designs. We apply this approach to the nested minimum end-compliance topology optimization using both the density-based method (non-convex) and the ground structure method (convex) by minimizing a weighted sum/average of the compliance over many load cases. Minimizing the weighted average of the compliance over many load cases requires the solution of the state equations for each load case in every optimization step to compute the (weighted) compliance and its gradient. Since the objective function and its gradient can be defined as the traces of symmetric matrices, we use the Hutchinson trace estimator, which provides the lowest variance, combined with the sample average approximation technique, to estimate both. This reduces the computational cost from $m\times N_\text{step}$ to roughly $n\times N_\text{step}$, where $m$ is the number of load cases, and $n \ll m$ is the sample size. We further propose a damping scheme for the stochastic algorithm, derived from simulated annealing to obtain fast convergence. We discuss the algorithmic parameters for our scheme and provide some information on how to choose these. 
 
The results for several generic examples and practical engineering designs demonstrate that the proposed stochastic algorithm provides high-quality designs at a drastically reduced computational cost. Based on the limited number and size of examples investigated, we show that the proposed stochastic algorithm may reduce the computational cost by a factor of up to 45, while obtaining a final compliance very close to that obtained by the standard algorithm, and occasionally even better. Our proposed damping scheme leads to fast convergence of the optimization. In addition, the combination of the discrete filter with the stochastic algorithm leads fewer linear solves with smaller systems, resulting in great computational efficiency.

The proposed stochastic algorithm is flexible and can be combined with any gradient-based update scheme. In this paper, we combine this technique with the optimality criteria, but other optimization methods can be used, such as the method of moving asymptotes (MMA). There are several important directions for future research. Although we have provided an important proof-of-concept, many questions remain open. One important question concerns optimal choices of parameters for both overall computational cost and quality of design. Future work should also consider dynamically varying the parameters. For example, if necessary, the quality of designs may be improved by increasing the sample size ($n$) as the minimum is approached. Another important question concerns the choice among stochastic optimization methods/approaches, both in terms of the overall methods as well as the choices in estimates and approximations. We need to further analyze what topology optimization problems can be solved efficiently using this approach. Finally, it would be useful to prove convergence to either a local or global minimum of the topology optimization problem under appropriate conditions. Since the weighted average of the compliance can be associated with uncertainties in the loads, this technique may be useful for design under uncertainty. Thus, we hope that the present paper will motivate the further use of stochastic sampling in various areas connected to (large scale) topology optimization.

%% file: Main_File-Stochastic_Sampling.bbl
\begin{thebibliography}{10}
\expandafter\ifx\csname url\endcsname\relax
  \def\url#1{\texttt{#1}}\fi
\expandafter\ifx\csname urlprefix\endcsname\relax\def\urlprefix{URL }\fi
\expandafter\ifx\csname href\endcsname\relax
  \def\href#1#2{#2} \def\path#1{#1}\fi

\bibitem{diaz1992}
A.~R. Diaz, M.~P. Bends\o{}e, Shape optimization of structures for multiple
  loading conditions using a homogenization method, Structural Optimization
  4~(1) (1992) 17--22.

\bibitem{bendsoe1994}
M.~P. Bends\o{}e, A.~Ben-Tal, J.~Zowe, Optimization methods for truss geometry
  and topology design, Structural optimization 7~(3) (1994) 141--159.

\bibitem{bental1997}
A.~Ben-Tal, A.~Nemirovski, Robust truss topology design via semidefinite
  programming, SIAM Journal on Optimization 7~(4) (1997) 991--1016.

\bibitem{bendsoe2003}
M.~P. Bends\o{}e, O.~Sigmund, Topology optimization: theory, methods, and
  applications, Springer, Berlin, Germany, 2003.

\bibitem{christiansen2001}
S.~Christiansen, P.~M., L.~Wynter, Stochastic bilevel programming in structural
  optimization, Structural and Multidisciplinary Optimization 21~(5) (2001)
  361--371.

\bibitem{achtziger1998}
W.~Achtziger, Multiple-load truss topology and sizing optimization: Some
  properties of minimax compliance, Journal of optimization theory and
  applications 98~(2) (1998) 255--280.

\bibitem{haber2012}
E.~Haber, M.~Chung, F.~Herrmann, An effective method for parameter estimation
  with {PDE} constraints with multiple right-hand sides, SIAM Journal on
  Optimization 22~(3) (2012) 739--757.

\bibitem{avron2010}
H.~Avron, P.~Maymounkov, S.~Toledo, Blendenpik: Supercharging {LAPACK}'s
  least-squares solver, SIAM Journal on Scientific Computing 32~(3) (2010)
  1217--1236.

\bibitem{Bourdin2001}
B.~Bourdin, Filters in topology optimization, International Journal for
  Numerical Methods in Engineering 50~(9) (2001) 2143--2158.

\bibitem{bensoe1989}
M.~P. Bends{\o}e, Optimal shape design as a material distribution problem,
  Structural optimization 1~(4) (1989) 193--202.

\bibitem{zhou1991}
M.~Zhou, G.~Rozvany, The {COC} algorithm, {P}art {II}: Topological, geometrical
  and generalized shape optimization, Computer Methods in Applied Mechanics and
  Engineering 89~(1-3) (1991) 309--336.

\bibitem{stolpe2001}
M.~Stolpe, K.~Svanberg, An alternative interpolation scheme for minimum
  compliance topology optimization, Structural and Multidisciplinary
  Optimization 22~(2) (2001) 116--124.

\bibitem{petersson1999}
J.~Petersson, A finite element analysis of optimal variable thickness sheets,
  SIAM Journal on Numerical Analysis 36~(6) (1999) 1759--1778.

\bibitem{Svanberg1984}
K.~Svanberg, On local and global minima in structural optimization, in:
  E.~Gallhager, R.~H. Ragsdell, O.~C. Zienkiewicz (Eds.), New Directions in
  Optimum Structural Design, John Wiley and Sons, Chichester, 1984, pp.
  327--341.

\bibitem{Achtziger1997}
W.~Achtziger, Topology optimization of discrete structures, Vol. 374 of
  International Centre for Mechanical Sciences, Springer Vienna, 1997, Ch.
  Topology Optimization in Structural Mechanics, pp. 57--100.

\bibitem{groenwold2008}
A.~A. Groenwold, L.~Etman, On the equivalence of optimality criterion and
  sequential approximate optimization methods in the classical topology layout
  problem, International journal for numerical methods in engineering 73~(3)
  (2008) 297--316.

\bibitem{hutchinson1989}
M.~Hutchinson, A stochastic estimator of the trace of the influence matrix for
  {L}aplacian smoothing splines, Communication in Statistics-Simulation and
  Computation 18~(3) (1989) 1059--1076.

\bibitem{avron2011}
H.~Avron, S.~Toledo, Randomized algorithms for estimating the trace of an
  implicit symmetric positive semi-definite matrix, Journal of the ACM 58~(2),
  article 8.

\bibitem{roosta2015}
F.~Roosta-Khorasani, U.~Ascher, Improved bounds on sample size for implicit
  matrix trace estimators, Foundations of Computational Mathematics 15~(5)
  (2015) 1187--1212.

\bibitem{shapiro2009}
A.~Shapiro, D.~Dentcheva, A.~Ruszczy\'{n}ski, Lectures on stochastic
  programming: modeling and theory, MPS-SIAM Series on Optimization, SIAM,
  Philadelphia, 2009.

\bibitem{anton2001}
A.~J. Kleywegt, A.~Shapiro, T.~H. de~Mello, The sample average approximation
  method for stochastic discrete optimization, SIAM Journal on Optimization
  12~(2) (2001) 479--502.

\bibitem{robbins1951}
H.~Robbins, S.~Monro, A stochastic approximation method, The Annals of
  Mathematical Statistics 22~(3) (1951) 400--407.

\bibitem{schraudolph1999}
N.~N. Schraudolph, Local gain adaptation in stochastic gradient descent, in:
  Intl. Conf. Artificial Neural Networks, Vol.~2, IEE, Edinburgh, Scotland,
  1999, pp. 569--574.

\bibitem{vishwanathan2006}
S.~V.~N. Vishwanathan, N.~N. Schraudolph, M.~W. Schmidt, K.~P. Murphy,
  Accelerated training of conditional random fields with stochastic gradient
  methods, in: Proceedings of the 23rd international conference on Machine
  learning, ACM, Pittsburgh, PA, 2006, pp. 969--976.

\bibitem{schraudolph2007}
N.~Schraudolph, J.~Yu, S.~G\"{u}nter, A stochastic quasi-newton method for
  online convex optimization, in: 11th Intl. Conf. on Artificial Intelligence
  and Statistics (AIstats), Soc. for Artificial Intelligence and Statistics,
  2007, pp. 433--440.

\bibitem{nemirovski2009}
A.~Nemirovski, A.~Juditsky, G.~Lan, A.~Shapiro, Robust stochastic approximation
  approach to stochastic programming, SIAM Journal on Optimization 19~(4)
  (2009) 1574--1609.

\bibitem{kirkpatrick1983}
S.~Kirkpatrick, C.~D. Gelatt~Jr, M.~P. Vecchi, Optimization by simulated
  annealing, Science 220~(4598) (1983) 671--680.

\bibitem{salamon2002}
P.~Salamon, P.~Sibani, R.~Frost, Facts, conjectures, and improvements for
  simulated annealing, SIAM Monographs on Mathematical Modeling and
  Computation, SIAM, Philadelphia, 2002.

\bibitem{svanberg1987}
K.~Svanberg, The method of moving asymptotes- a new method for structural
  optimization, International Journal for Numerical Methods in Engineering
  24~(2) (1987) 359--373.

\bibitem{Ramos2016}
A.~S. Ramos~Jr, G.~H. Paulino, {Filtering structures out of ground structures}
  -- {a} discrete filtering tool for structural design optimization, Structural
  and Multidisciplinary Optimization 54 (2016) 95--116.

\bibitem{Talischi2012a}
C.~Talischi, G.~H. Paulino, A.~Pereira, I.~F.~M. Menezes, {PolyTop: {A} Matlab
  implementation of a general topology optimization framework using
  unstructured polygonal finite element meshes}, Structural and
  Multidisciplinary Optimization 45~(3) (2012) 329--357.
\newblock \href {http://dx.doi.org/10.1007/s00158-011-0696-x}
  {\path{doi:10.1007/s00158-011-0696-x}}.

\bibitem{liu2014}
K.~Liu, A.~Tovar, An efficient {3D} topology optimization code written in
  {MATLAB}, Structural and Multidisciplinary Optimization 50~(6) (2014)
  1175--1196.

\bibitem{zegard2016}
T.~Zegard, G.~H. Paulino, Bridging topology optimization and additive
  manufacturing, Structural and Multidisciplinary Optimization 53~(1) (2016)
  175--192.

\bibitem{zegard2014}
T.~Zegard, G.~H. Paulino, {GRAND} -- {G}round structure based topology
  optimization for arbitrary {2D} domains using {MATLAB}, Structural and
  Multidisciplinary Optimization 50~(5) (2014) 861--882.

\bibitem{zegard2015a}
T.~Zegard, G.~H. Paulino, {GRAND3} -- {G}round structure based topology
  optimization for arbitrary {3D} domains using {MATLAB}, Structural and
  Multidisciplinary Optimization 52~(6) (2015) 1161--1184.

\bibitem{christensen2009}
P.~Christensen, A.~Klarbring, An introduction to structural optimization, Vol.
  153, Springer Science \& Business Media, 2009.

\bibitem{zhang2016a}
X.~Zhang, S.~Maheshwari, A.~S. Ramos~Jr, G.~H. Paulino, Macroelement and
  macropatch approaches to structural topology optimization using the ground
  structure method, ASCE Journal of Structural Engineering.\href
  {http://dx.doi.org/10.1061/(ASCE) ST.1943-541X.0001524.}
  {\path{doi:10.1061/(ASCE) ST.1943-541X.0001524.}}

\bibitem{zhang2016}
X.~Zhang, A.~S. Ramos~Jr, G.~H. Paulino, A discrete filter scheme for material
  nonlinear topology design using the ground structure method, Structural and
  Multidisciplinary Optimization. Submitted.

\bibitem{canyantower}
{Cayan Tower (Infinity Tower)}, {URL}: http://cayan.net.

\end{thebibliography}
